\documentclass[aap,mtplusscr,seceqn,secthm,dvips]{arximspdf}
\usepackage{mathrsfs}

\doi{10.1214/105051606000000871}
\volume{17}
\issue{2}
\pubyear{2007}
\firstpage{688}
\lastpage{715}

\makeatletter

\newproclaim{notation}{Notation}
\newtheorem{lem}[thm]{Lemma}
\newtheorem{cor}[thm]{Corollary}
\newproclaim{rem}[thm]{Remark}
\newtheorem{prop}[thm]{Proposition}
\newcommand{\mall}{\boldsymbol{\sf D}}
\newcommand{\nab}{\boldsymbol{\nabla}}

\makeatother

\begin{document}
\begin{frontmatter}

\title{Model robustness of finite state nonlinear filtering over the
infinite time horizon}
\runtitle{Model robustness of finite state nonlinear filtering}

\begin{aug}
\author[A]{\fnms{Pavel} \snm{Chigansky}\ead[label=e1]{pavel.chigansky@weizmann.ac.il}\thanksref{f1}}
and
\author[B]{\fnms{Ramon} \snm{van Handel}\corref{}\ead[label=e2]{ramon@its.caltech.edu}\thanksref{f2}}
\runauthor{P. Chigansky and R. van Handel}
\affiliation{Weizmann Institute of Science and California Institute of Technology}
\address[A]{Department of Mathematics \\
The Weizmann Institute of Science \\
Rehovot 76100 \\
Israel \\ \printead{e1}} 
\address[B]{Physical Measurement and Control 266-33 \\
California Institute of Technology \\
Pasadena, California 91125 \\
USA \\ \printead{e2}}
\thankstext{f1}{Supported by a grant from the Israeli Science Foundation.}
\thankstext{f2}{Supported by the ARO under Grants DAAD19-03-1-0073 and W911NF-06-1-0378.}
\end{aug}

\received{\smonth{4} \syear{2006}}
\revised{\smonth{11} \syear{2006}}

\begin{abstract}
We investigate the robustness of nonlinear filtering for continuous time
finite state Markov chains, observed in white noise, with respect to
misspecification of the model parameters. It is shown that the distance
between the optimal filter and that with incorrect model parameters
converges to zero uniformly over the infinite time interval as the
misspecified model converges to the true model, provided the signal obeys
a mixing condition.  The filtering error is controlled through the
exponential decay of the derivative of the nonlinear filter with respect
to its initial condition.  We allow simultaneously for misspecification of
the initial condition, of the transition rates of the signal, and of the
observation function.  The first two cases are treated by relatively
elementary means, while the latter case requires the use of Skorokhod
integrals and tools of anticipative stochastic calculus.
\end{abstract}



\end{frontmatter}

\section{Introduction}
\label{sec:intro}

The theory of nonlinear filtering concerns the estimation of a signal
corrupted by white noise, and has diverse applications in target tracking,
signal processing, automatic control, finance, and so on.  The basic
setting of the theory involves a Markov signal process, for example, the solution
of a (nonlinear) stochastic differential equation or a finite-state Markov
process, observed in independent corrupting noise.  The calculation of the
resulting filters is a classical topic in stochastic analysis
\cite{liptser}.  Of course, the filtering equations will depend explicitly
on the model chosen for the signal process and observations; in almost all
realistic applications, however, the model that underlies the filter is
only an approximation of the true system that generates the observations.
In order for the theory to be practically useful, it is important to
establish that the filtered estimates are not too sensitive to the choice
of underlying model.

Continuity with respect to the model parameters of nonlinear filtering
estimates on a fixed \textit{finite} time interval is well established,
for example, \cite{BhatKalKar1,BhatKalKar2,GuoYin}; generally speaking, it is known
that the error incurred in a finite time interval due to the choice of
incorrect model parameters can be made arbitrarily small if the model
parameters are chosen sufficiently close to those of the true model.  As
the corresponding error bounds grow rapidly with the length of the time
interval, however, such estimates are of little use if we are interested
in robustness of the filter over a long period of time.  One would like to
show that the approximation errors do not accumulate, so that the error
remains bounded uniformly over an \textit{infinite} time interval.

The model robustness of nonlinear filters on the infinite time horizon was
investigated in discrete time in \cite{kushner1,legland1,legland2,papa}.
The key idea that allows one to control the accumulation of
approximation errors is the asymptotic stability property of many
nonlinear filters, which is the focus of much recent work (see
\cite{atar,BaxChiLip} and the references therein) and can be
summarized as follows.  The optimal nonlinear filter is a recursive
equation that is initialized with the true distribution of the signal
process at the initial time.  If the filter is initialized with a
different distribution, then the resulting filtered estimates are no
longer optimal (in the least-squares sense).  The filter is called
asymptotically stable if the solution of the wrongly initialized filter
converges to the solution of the correctly initialized filter at large
times; that is, the filter ``forgets'' its initial condition after a period
of observation.

Using an approximate filter rather than the optimal filter is equivalent
to using the optimal filter where we make an approximation error after
every time step. Now suppose the optimal filter forgets its initial
condition at an exponential rate; then also the approximation error at
each time step is forgotten at an exponential rate, and the errors cannot
accumulate in time.  If the approximation error at each time step is
bounded (finite time robustness), then the total approximation error will
be bounded uniformly in time.  Model robustness on the infinite time
horizon is thus a consequence of finite time robustness together with the
exponential forgetting property of the filter. This is precisely the
method used in \cite{kushner1,legland1,legland2,papa}, and its
implementation is fairly straightforward once bounds on the exponential
forgetting rate of the filter have been obtained.  However, the method
used there does not extend to nonlinear filtering in continuous time; even
the continuous time model with point process observations studied in
\cite{kushner1}, though more involved, reduces essentially to discrete
(but random) observation times.  The continuous time case requires
different tools, which we develop in this paper in the setting of
nonlinear filtering of a finite-state Markov signal process. (We also
mention \cite{kushner2}, where a different but related problem is solved.)

We consider the following filtering setup.  The signal process
$X=(X_t)_{t\ge 0}$ is a continuous time, homogeneous Markov chain with
values in the finite alphabet $\mathbb{S}=\{a_1,\ldots,a_d\}$, transition
intensities matrix $\Lambda=(\lambda_{ij})$ and initial distribution
$\nu^i=\mathbf{P}(X_0=a_i)$.  The observation process $Y=(Y_t)_{t\ge 0}$
is given by
\begin{equation}
\label{eq:observation}
Y_t=\int_0^t \mathrm{h}(X_s)\,ds+B_t,
\end{equation}
where $\mathrm{h}\dvtx \mathbb{S}\to\mathbb{R}$ is the observation function
[we will also write $h^i=\mathrm{h}(a_i)$] and $B$ is a Wiener process that is
independent of $X$. The filtering problem for this model concerns the
calculation of the conditional probabilities
$\pi_t^i=\mathbf{P}(X_t=a_i|\mathscr{F}_t^Y)$ from the observations
$\{Y_s\dvtx s\le t\}$, where $\mathscr{F}_t^Y=\sigma\{Y_s\dvtx s\le t\}$.  It is
well known that $\pi_t$ satisfies the Wonham equation \cite{liptser,wonham}
\begin{equation}
\label{eq:wonham}
d\pi_t=\Lambda^*\pi_t\,dt+(H-h^*\pi_t)\pi_t\,(dY_t-h^*\pi_t\,dt),
\qquad \pi_0=\nu,
\end{equation}
where $x^*$ denotes the transpose of $x$ and $H={\rm diag}\,h$.  Note that
the Wonham equation is initialized with the true distribution of $X_0$; we
will denote by $\pi_t(\mu)$ the solution of the Wonham equation at time
$t$ with an arbitrary initial distribution $\pi_0=\mu$, and by
$\pi_{s,t}(\mu)$ the solution of the Wonham equation at time $t\ge s$ with
the initial condition $\pi_s=\mu$.  In \cite{atar,BaxChiLip} the
exponential forgetting property of the Wonham filter was established as
follows: the $\ell_1$-distance $|\pi_t(\mu)-\pi_t(\nu)|$ decays
exponentially a.s., provided the initial distributions are equivalent
$\mu\sim\nu$ and that the {\em mixing condition} $\lambda_{ij}>0$ $\forall
i\ne j$ is satisfied.  Now consider the Wonham filter with incorrect model
parameters:
\begin{equation}
\label{eq:wrongwonham}
d\tilde\pi_t=\tilde\Lambda^*\tilde\pi_t\,dt+(\tilde H-\tilde h^*\tilde\pi_t)\tilde\pi_t
(dY_t-\tilde h^*\tilde\pi_t\,dt), \qquad \tilde\pi_0=\nu,
\end{equation}
where $\tilde\Lambda$ and $\tilde h$ denote a transition intensities matrix and
observation function that do not match the underlying signal-observation
model $(X,Y)$, $\tilde H=\operatorname{diag}\tilde h$, and we denote by $\tilde\pi_t(\mu)$ the
solution of this equation with initial condition $\tilde\pi_0=\mu$ and by
$\tilde\pi_{s,t}(\mu)$ the solution with $\tilde\pi_s=\mu$.  The following is the
main result of this paper.

\begin{thm}
\label{thm:mainresult}
Suppose $\nu^i,\mu^i>0$ $\forall i$ and
$\lambda_{ij},\tilde\lambda_{ij}>0$ $\forall i\ne j$.  Then
\[
\sup_{t\ge 0}\mathbf{E}\|\tilde\pi_t(\mu)-\pi_t(\nu)\|^2
\le C_1\,|\mu-\nu|+ C_2\,|\tilde h-h|+ C_3|\tilde\Lambda^*-\Lambda^*|,
\]
where $|\tilde\Lambda^*-\Lambda^*|=\sup\{|(\tilde\Lambda^*-\Lambda^*)\tau|\dvtx\tau^i>0\ \forall
i,|\tau|=1\}$ and the quantities $C_1,C_2,C_3$ are bounded on any compact
subset of parameters $\{(\nu,\Lambda,h,\mu,\tilde\Lambda,\tilde h):
\nu^i,\mu^i>0\ \forall i,|\nu|=|\mu|=1,\lambda_{ij},\tilde\lambda_{ij}>0~
\forall i\ne j,\sum_j\lambda_{ij}=\sum_j\tilde\lambda_{ij}=0~\forall i\}$.
Additionally we have the asymptotic estimate
\[
\limsup_{t\to\infty}\mathbf{E}\|\tilde\pi_t(\mu)-\pi_t(\nu)\|^2
\le C_2|\tilde h-h|+ C_3|\tilde\Lambda^*-\Lambda^*|.
\]
In particular, this implies that if $\nu^i>0$ $\forall i$,
$\lambda_{ij}>0$ $\forall i\ne j$, then
\[
\lim_{(\tilde h,\tilde\Lambda,\mu)\to(h,\Lambda,\nu)}
\sup_{t\ge 0}\mathbf{E}\|\tilde\pi_t(\mu)-\pi_t(\nu)\|=
\lim_{(\tilde h,\tilde\Lambda)\to(h,\Lambda)}
\limsup_{t\to\infty}\mathbf{E}\|\tilde\pi_t(\mu)-\pi_t(\nu)\|=0.
\]\
\end{thm}

Let us sketch the basic idea of the proof.  Rather than considering
the Wonham filter, let us demonstrate the idea using the following simple
caricature of a filtering equation.  Consider a smooth ``observation''
$y_t$ and a ``filter'' whose state $x_t$ is propagated by the
ordinary differential equation $dx_t/dt=f(x_t,y_t)$.  Similarly, we
consider the ``approximate filter'' $d\tilde x_t/dt=\tilde f(\tilde
x_t,y_t)$ and assume that everything is sufficiently smooth, so that for
fixed $y$ both equations generate a two-parameter flow
$x_t=\varphi_{s,t}^y(x_s)$, $\tilde x_t=\tilde\varphi_{s,t}^y(\tilde x_s)$.
The following calculation is straightforward:
\begin{eqnarray*}
\varphi^y_{0,t}(x)-\tilde\varphi^y_{0,t}(x) &=&
\int_0^t\frac{d}{ds}[\tilde\varphi^y_{s,t}(\varphi^y_{0,s}(x))]\,ds\\
\\
&=&\int_0^t D\tilde\varphi^y_{s,t}(\varphi^y_{0,s}(x))\cdot
\bigl(f(\varphi^y_{0,s}(x),y_s)-\tilde f(\varphi^y_{0,s}(x),y_s)\bigr)\,ds,
\end{eqnarray*}
where $D\tilde\varphi^y_{s,t}(x)\cdot v$ denotes the directional
derivative of $\tilde\varphi^y_{s,t}(x)$ in the direction $v$.
Hence we obtain the following estimate on the approximation error:
\[
|\varphi^y_{0,t}(x)-\tilde\varphi^y_{0,t}(x)|\le\int_0^t
|D\tilde\varphi^y_{s,t}(\varphi^y_{0,s}(x))|~
|f(\varphi^y_{0,s}(x),y_s)-\tilde f(\varphi^y_{0,s}(x),y_s)|\,ds.
\]
Now suppose that $|f(\cdot,\cdot)-\tilde f(\cdot,\cdot)|\le K$, where
$K\to 0$ as $\tilde f\to f$; this is an expression of finite-time
robustness, as it ensures that
$|\varphi^y_{0,t}(x)-\tilde\varphi^y_{0,t}(x)|\le Kt\to 0$ (for fixed $t$)
as $\tilde f\to f$.  Suppose furthermore that we can establish a
bound of the form $|D\tilde\varphi^y_{s,t}(\cdot)|\le Ce^{-\lambda(t-s)}$,
that is an infinitesimal perturbation to the initial condition is forgotten
at an exponential rate.  Then the estimate above is uniformly bounded and
converges to zero uniformly in time as $\tilde f\to f$.  Conceptually this
is similar to the logic used in discrete time, but we have to replace
the exponential forgetting of the initial condition by the requirement
that the derivative of the filter with respect to its initial condition
decays exponentially.

Returning to the Wonham filter, this procedure can be implemented in a
fairly straightforward way if $\tilde h=h$.  In this case, most of the work
involves finding a suitable estimate on the exponential decay of the
derivative of the filter with respect to its initial condition; despite
the large number of results on filter stability, such estimates are not
available in the literature to date.  We obtain estimates by adapting
methods from \cite{BaxChiLip}, together with uniform estimates of the
concentration of the optimal filter near the boundary of the simplex.

The general case with $\tilde h\ne h$ is significantly more involved. The
problem is already visible in the simple demonstration above.  Note that
the integrand on the right-hand side of the error estimate is not adapted;
it depends on the observations on the entire interval $[0,t]$.  As the
Wonham filter is defined in terms of an It\^o-type stochastic integral,
this will certainly get us into trouble.  When $\tilde h=h$ the stochastic
integral cancels in the error bound and the problems are kept to a
minimum; in the general case, however, we are in no such luck.
Nonetheless this problem is not prohibitive, but it requires us to
use the stochastic calculus for anticipating integrands developed by
Nualart and Pardoux \cite{nualartpardoux,nualart} using Skorokhod
integrals rather than It\^o integrals and using Malliavin calculus tools.

An entirely different application of the Malliavin calculus to
problems of filter stability can be found in \cite{daprato}.

The remainder of this paper is organized as follows.  In Section
\ref{sec:prelims} we prove some regularity properties of the solution of
the Wonham equation.  We also demonstrate the error estimate discussed
above in the simpler case $\tilde h=h$, and comment on the more general
applicability of such a bound.  In Section \ref{sec:forgetting} we obtain
exponential bounds on the derivative of the Wonham filter with respect to
its initial condition.  Section \ref{sec:robustness} treats the general
case $\tilde h\ne h$ using anticipative stochastic calculus; some of the
technical estimates appear in Appendix \ref{sec:proofs}.  Finally,
Appendix \ref{sec:malliavin} contains a brief review of the results from
the Malliavin calculus and anticipative stochastic calculus that are
needed in the proofs.
\begin{notation*}
The signal-observation pair $(X,Y)$ is defined on the standard probability
space $(\Omega,\mathscr{F},\mathbf{P})$.  The expectation with respect to
$\mathbf{P}$ is denoted by $\mathbf{E}$ or sometimes $\mathbf{E_P}$.  For
$x\in\mathbb{R}^d$, we denote by $|x|$ the $\ell_1$-norm, by
$\|x\|$ the $\ell_2$-norm, and by $\|x\|_p$ the $\ell_p$-norm.  We write
$x\succ y$ (resp. $\prec,\succeq,\preceq$) if $x_i>y_i$ ($<,\ge,\le$) $\forall i$.

The following spaces will be used throughout.  Probability distributions
on $\mathbb{S}$ are elements of the simplex
$\Delta^{d-1}=\{x\in\mathbb{R}^d\dvtx x\succeq 0, |x|=1\}$. Usually, we will be
interested in the interior of the simplex
$\mathcal{S}^{d-1}=\{x\in\mathbb{R}^d\dvtx x\succ 0, |x|=1\}$.  The space of
vectors tangent to $\mathcal{S}^{d-1}$ is denoted by
$T\mathcal{S}^{d-1}=\{x\in\mathbb{R}^d\dvtx \sum_ix_i=0\}$.
Finally, we will denote the positive orthant by
$\mathbb{R}^d_{++}=\{x\in\mathbb{R}^d\dvtx x\succ 0\}$.
\end{notation*}

\section{Preliminaries}
\label{sec:prelims}

Equation (\ref{eq:wonham}) is a nonlinear equation for the conditional
distribution $\pi_t$.  It is well known however (e.g. \cite{elliott})
that $\pi_t$ can also be calculated in a linear fashion:
$\pi_t=\rho_t/|\rho_t|$, where the unnormalized density $\rho_t$ is
propagated by the Zakai equation
\begin{equation}
d\rho_t = \Lambda^*\rho_t\,dt + H\rho_t\,dY_t,\qquad \rho_0=\nu.
\end{equation}
We will repeatedly exploit this representation in what follows.  As
before $\rho_t(\mu)$ and $\rho_{s,t}(\mu)$ ($t\ge s$) denote the solution
of the Zakai equation at time $t$ with the initial condition $\rho_0=\mu$
and $\rho_s=\mu$, respectively, and $\pi_{s,t}(\mu)=
\rho_{s,t}(\mu)/|\rho_{s,t}(\mu)|$.

We also recall the following interpretation of the norm $|\rho_t|$ of the
unnormalized conditional distribution.  If we define a new measure
$\mathbf{Q}\sim\mathbf{P}$ through
\begin{equation}
\label{eq:refmeasure}
\frac{d{\bf P}}{d{\bf Q}}=|\rho_t(\nu)|=|\rho_t|,
\end{equation}
then under $\mathbf{Q}$ the observation process $Y_t$ is an
$\mathscr{F}_t^Y$-Wiener process.  This observation will be used in
Section \ref{sec:robustness} to apply the Malliavin calculus.

The main goal of this section is to establish some regularity properties
of the solutions of the Wonham and Zakai equations.  In particular, as we
will want to calculate the derivative of the filter with respect to its
initial condition, we have to establish that $\pi_{s,t}(\mu)$ is in fact
differentiable.  We will avoid problems at the boundary of the simplex by
disposing of it alltogether: we begin by proving that if
$\mu\in\mathcal{S}^{d-1}$, then a.s. $\pi_{s,t}(\mu)\in\mathcal{S}^{d-1}$
for all times $t>s$.
\begin{lem}
${\bf P}(\rho_{s,t}(\mu)\in\mathbb{R}^d_{++}
\mbox{ for all }\mu\in\mathbb{R}^d_{++},~0\le s\le t<\infty)=1$.
\end{lem}

\begin{pf}
The following variant on the pathwise filtering method reduces the Zakai
equation to a random differential equation.  First, we write
$\Lambda^*=S+T$ where $S$ is the diagonal matrix with
$S_{ii}=\lambda_{ii}$.  Note that the matrix $T$ has only nonnegative
entries.  We now perform the transformation $f_{s,t}(\mu)=L_{s,t}\rho_{s,t}(\mu)$ where
\[
L_{s,t}=\exp\bigl(\bigl(\tfrac{1}{2}H^2-S\bigr)(t-s)-H(Y_t-Y_s)\bigr).
\]
Then $f_{s,t}(\mu)$ satisfies
\begin{equation}
\label{eq:robustfilt}
\frac{df_{s,t}}{dt}=L_{s,t}TL_{s,t}^{-1}\,f_{s,t},\qquad f_{s,s}=\mu.
\end{equation}
Let $\Omega_c\subset\Omega$, ${\bf P}(\Omega_c)=1$ be a set such that
$t\mapsto B_t(\omega)$ is continuous for every $\omega\in\Omega_c$. Then
$t\mapsto L_{s,t}$, $t\mapsto L_{s,t}^{-1}$ are continuous in $t$ and have
strictly positive diagonal elements for every $\omega\in\Omega_c$. By
standard arguments, there exists for every $\omega\in\Omega_c$,
$\mu\in\mathbb{R}^d$ and $s\ge 0$ a unique solution $f_{s,t}(\mu)$ to equation
(\ref{eq:robustfilt}) where $t\mapsto f_{s,t}(\mu)$ is a $C^1$-curve.
Moreover, note that $L_{s,t}TL_{s,t}^{-1}$ has nonnegative matrix elements
for every $\omega\in\Omega_c$, $s\le t<\infty$.  Hence if
$\mu\in\mathbb{R}^d_{++}$ then clearly $f_{s,t}(\mu)$ must be
nondecreasing, that is, $f_{s,t}\succeq f_{s,r}$ for every $t\ge r\ge s$ and
$\omega\in\Omega_c$. But then $\mathbb{R}^d_{++}$ must be forward
invariant under equation (\ref{eq:robustfilt}) for every $\omega\in\Omega_c$,
and as $L_{s,t}$ has strictly positive diagonal elements the result follows.
\end{pf}

\begin{cor}
${\bf P}(\pi_{s,t}(\mu)\in\mathcal{S}^{d-1}
\mbox{ for all }\mu\in\mathcal{S}^{d-1},~0\le s\le t<\infty)=1$.
\end{cor}

Let us now investigate the map $\rho_{s,t}(\mu)$.  As this map is linear
in $\mu$, we can write $\rho_{s,t}(\mu)=U_{s,t}\mu$ a.s.\ where the
$d\times d$ matrix $U_{s,t}$ is the solution of
\begin{equation}
\label{eq:zakaiprop}
dU_{s,t} = \Lambda^*U_{s,t}\,dt+ HU_{s,t}\,dY_t, \qquad U_{s,s}=I.
\end{equation}
The following lemma establishes that $U_{s,t}$ defines a linear stochastic
flow in $\mathbb{R}^d$.

\begin{lem}
\label{lem:zakaiflow}
For a.e. $\omega\in\Omega$ \emph{(i)} $\rho_{s,t}(\mu)=U_{s,t}\mu$ for
all $s\le t$; \emph{(ii)} $U_{s,t}$ is continuous in $(s,t)$;
\emph{(iii)} $U_{s,t}$ is invertible for all $s\le t$, where $U_{s,t}^{-1}$ is given by
\begin{equation}
\label{eq:zakaiinv}
dU_{s,t}^{-1} = -U_{s,t}^{-1}\Lambda^*\,dt +U_{s,t}^{-1}H^2\,dt
-U_{s,t}^{-1}H\,dY_t, \qquad U_{s,s}^{-1}=I;
\end{equation}
\emph{(iv)} $U_{r,t}U_{s,r}=U_{s,t}$ (and hence
$U_{s,t}U_{s,r}^{-1}=U_{r,t}$) for all $s\le r\le t$.
\end{lem}

\begin{pf}
Continuity of $U_{s,t}$ (and $U_{s,t}^{-1}$) is a standard property of
solution of Lipschitz stochastic differential equations.  Invertibility of
$U_{0,t}$ for all $0\le t<\infty$ is established in
\cite{protter}, page 326, and it is evident that
$U_{s,t}=U_{0,t}U_{0,s}^{-1}$ satisfies equation (\ref{eq:zakaiprop}).
The remaining statements follow, where we can use continuity to remove the
time dependence of the exceptional set as in the proof of
\cite{protter}, page 326.
\end{pf}

We now turn to the properties of the map $\pi_{s,t}(\mu)$.

\begin{lem}
The Wonham filter generates a smooth stochastic semiflow in
$\mathcal{S}^{d-1}$, that is, the solutions $\pi_{s,t}(\mu)$ satisfy the
following conditions:
\begin{enumerate}
\item For a.e. $\omega\in\Omega$, $\pi_{s,t}(\mu)=
\pi_{r,t}(\pi_{s,r}(\mu))$ for all $s\le r\le t$ and $\mu$.
\item For a.e. $\omega\in\Omega$, $\pi_{s,t}(\mu)$ is continuous
in $(s,t,\mu)$.
\item For a.e. $\omega\in\Omega$, the injective map $\pi_{s,t}(\cdot)\dvtx
\mathcal{S}^{d-1}\to\mathcal{S}^{d-1}$ is $C^\infty$ for all $s\le t$.
\end{enumerate}
\end{lem}

\begin{pf}
For $x\in\mathbb{R}^d_{++}$ define $\Sigma(x)=x/|x|$, so that
$\pi_{s,t}(\mu)=\Sigma(\rho_{s,t}(\mu))$ ($\mu\in\mathcal{S}^{d-1}$).
Note that $\Sigma$ is smooth on $\mathbb{R}^d_{++}$.  Hence continuity in
$(s,t,\mu)$ and smoothness with respect to $\mu$ follow directly from the
corresponding properties of $\rho_{s,t}(\mu)$.  The semiflow property
$\pi_{s,t}(\mu)=\pi_{r,t}(\pi_{s,r}(\mu))$ follows directly from Lemma
\ref{lem:zakaiflow}.  It remains to prove injectivity.

Suppose that $\pi_{s,t}(\mu)=\pi_{s,t}(\nu)$ for some
$\mu,\nu\in\mathcal{S}^{d-1}$. Then
$U_{s,t}\mu/|U_{s,t}\mu|=U_{s,t}\nu/|U_{s,t}\nu|$, and as $U_{s,t}$ is
invertible we have $\mu=(|U_{s,t}\mu|/|U_{s,t}\nu|)\nu$.  But as $\mu$ and
$\nu$ must lie in $\mathcal{S}^{d-1}$, it follows that $\mu=\nu$.  Hence
$\pi_t(\cdot)$ is injective.
\end{pf}

\begin{rem}
The results in this section hold identically if we replace $\Lambda$
by $\tilde\Lambda$, $h$ by $\tilde h$.  We will use the obvious notation
$\tilde\pi_{s,t}(\mu)$, $\tilde\rho_{s,t}(\mu)$, $\tilde U_{s,t}$, and so on.
\end{rem}

We finish this section by obtaining an expression for the approximation
error in the case $\tilde h=h$; in fact, we will demonstrate the bound for
this simple case in a more general setting than is considered in the
following.  Rather than considering the approximate Wonham filter with
modified $\Lambda$, consider the equation
\begin{equation}
\label{eq:filterapproximation}
d\breve\pi_t=f(\breve\pi_t)\,dt +(H-h^*\breve\pi_t)\breve\pi_t
(dY_t-h^*\breve\pi_t\,dt),
\qquad \breve\pi_0=\mu\in\mathcal{S}^{d-1},
\hspace*{-15pt}
\end{equation}
where $f\dvtx\mathcal{S}^{d-1}\to T\mathcal{S}^{d-1}$ is chosen in such a
way that this equation has a strong solution and $\inf\{t>0\dvtx
\breve\pi_t\notin\mathcal{S}^{d-1}\}=\infty$ a.s.  In the sequel
we consider the case $f(\breve\pi)=\tilde\Lambda\breve\pi$, which clearly satisfies
the requirements.  We formulate the more general result here, as it might
be of interest in other contexts (see Remark \ref{rem:projf}).

\begin{prop}
\label{pro:simplebound}
Let $\breve\pi_t$ be as above.  Then the difference between $\breve\pi_t$
and the Wonham filter started at $\mu$ is a.s.\ given by
\[
\breve\pi_t-\pi_t(\mu)=\int_0^t
D\pi_{s,t}(\breve\pi_s)\cdot\bigl(f(\breve\pi_s)-\Lambda^*\breve\pi_s\bigr)\,ds,
\]
where $D\pi_{s,t}(\mu)\cdot v$ is the derivative of
$\pi_{s,t}(\mu)$ in the direction $v\in T\mathcal{S}^{d-1}$.
\end{prop}

\begin{pf}
Define the (scalar) process $\Gamma_t$ by
\[
\Gamma_t=\exp\biggl(\int_0^th^*\breve\pi_s\,dY_s
 -\tfrac{1}{2}\int_0^t(h^*\breve\pi_s)^2\,ds\biggr).
\]
Using It\^o's rule, we evaluate
\begin{equation}
\label{eq:itodropsout}
\frac{d}{ds}(\Gamma_sU_{0,s}^{-1}\breve\pi_s)=
\Gamma_sU_{0,s}^{-1}\bigl(f(\breve\pi_s)-\Lambda^*\breve\pi_s\bigr).
\end{equation}
Multiplying both sides by $U_{0,t}$, we obtain
\[
\frac{d}{ds}(\Gamma_sU_{s,t}\breve\pi_s)=
\Gamma_sU_{s,t}\bigl(f(\breve\pi_s)-\Lambda^*\breve\pi_s\bigr).
\]
Now introduce as before the map
$\Sigma\dvtx\mathbb{R}^{d}_{++}\to\mathcal{S}^{d-1}$, $\Sigma(x)=x/|x|$, which
is smooth on $\mathbb{R}^{d}_{++}$.  Define the matrix $D\Sigma(x)$ with elements
\[
[D\Sigma(x)]^{ij}=\frac{\partial\Sigma^i(x)}{\partial x^j}=
\frac{1}{|x|}[\delta_{ij}-\Sigma^i(x)].
\]
Note that $\Sigma(\alpha x)=\Sigma(x)$ for any $\alpha>0$.  Hence
\[
\frac{d}{ds}\Sigma(U_{s,t}\breve\pi_s)=
\frac{d}{ds}\Sigma(\Gamma_sU_{s,t}\breve\pi_s)=
D\Sigma(\Gamma_sU_{s,t}\breve\pi_s)
\frac{d}{ds}(\Gamma_sU_{s,t}\breve\pi_s).
\]
But then we have, using $D\Sigma(\alpha x)=\alpha^{-1}D\Sigma(x)$ ($\alpha>0$),
\begin{eqnarray*}
\frac{d}{ds}\Sigma(U_{s,t}\breve\pi_s) &=&
D\Sigma(\Gamma_sU_{s,t}\breve\pi_s)
\Gamma_sU_{s,t}\bigl(f(\breve\pi_s)-\Lambda^*\breve\pi_s\bigr)
\\
&=& D\Sigma(U_{s,t}\breve\pi_s)
U_{s,t}\bigl(f(\breve\pi_s)-\Lambda^*\breve\pi_s\bigr).
\end{eqnarray*}
On the other hand, we obtain from the representation
$\pi_{s,t}(\mu)=\Sigma(U_{s,t}\mu)$
\[
D\pi_{s,t}(\mu)\cdot v= D\Sigma(U_{s,t}\mu)U_{s,t}v,
\qquad \mu\in\mathcal{S}^{d-1},~v\in T\mathcal{S}^{d-1}.
\]
Note that $f(\breve\pi_s)-\Lambda^*\breve\pi_s\in T\mathcal{S}^{d-1}$
as we required that $f\dvtx\mathcal{S}^{d-1}\to T\mathcal{S}^{d-1}$,
so that $D\Sigma(U_{s,t}\breve\pi_s)U_{s,t}(f(\breve\pi_s)-
\Lambda^*\breve\pi_s)=D\pi_{s,t}(\breve\pi_s)\cdot
(f(\breve\pi_s)-\Lambda^*\breve\pi_s)$.  Finally, note that
\[
\int_0^t\frac{d}{ds}\Sigma(U_{s,t}\breve\pi_s)\,ds=
\Sigma(\breve\pi_t)-\Sigma(U_{0,t}\breve\pi_0)=\breve\pi_t-\pi_t(\mu),
\]
and the proof is complete.
\end{pf}

\begin{cor}
\label{cor:filterapproximation}
The following estimate holds:
\[
|\breve\pi_t-\pi_t(\mu)|\le\int_0^t
|D\pi_{s,t}(\breve\pi_s)| |f(\breve\pi_s)-\Lambda^*\breve\pi_s|\,ds,
\]
where $|D\pi_{s,t}(\mu)|=\sup\{|D\pi_{s,t}(\mu)\cdot v|\dvtx
v\in T\mathcal{S}^{d-1},~|v|=1\}$.  Moreover
\[
|\breve\pi_t-\pi_t(\nu)|\le|\pi_t(\mu)-\pi_t(\nu)|+
\int_0^t |D\pi_{s,t}(\breve\pi_s)| |f(\breve\pi_s)-\Lambda^*\breve\pi_s| \,ds.
\]
\end{cor}

\begin{rem}
\label{rem:projf}
Corollary \ref{cor:filterapproximation} suggests that the method used here
could be applicable to a wider class of filter approximations than those
obtained by misspecification of the underlying model.  In particular, in
the infinite-dimensional setting it is known \cite{brigo} that by
projecting the filter onto a properly chosen finite-dimensional manifold,
one can obtain finite-dimensional approximate filters that take a form
very similar to equation (\ref{eq:filterapproximation}).  In order to obtain
useful error bounds for such approximations one would need to have a
fairly tight estimate on the derivative of the filter with respect to its
initial condition.  Unfortunately, worst-case estimates of the type
developed in Section \ref{sec:forgetting} are not sufficiently tight to
give quantitative results on the approximation error, even in the
finite-state case.  In the remainder of the article we will restrict
ourselves to studying the robustness problem.
\end{rem}

In the following, it will be convenient to turn around the role of the
exact and approximate filters in Corollary \ref{cor:filterapproximation},
that is, we will use the estimate
\begin{equation}
\label{eq:simplebound}
|\pi_t(\nu)-\tilde\pi_t(\mu)|\le |\tilde\pi_t(\nu)-\tilde\pi_t(\mu)|+ \int_0^t
|D\tilde\pi_{s,t}(\pi_s)| |(\Lambda^*-\tilde\Lambda^*)\pi_s|\,ds,
\hspace*{-20pt}
\end{equation}
which holds provided $\tilde h=h$.  The proof is identical to the one given above.

\section{Exponential estimates for the derivative of the filter}
\label{sec:forgetting}

In order for the bound equation (\ref{eq:simplebound}) to be useful, we
must have an exponential estimate for $|D\tilde\pi_{s,t}(\cdot)|$.  The goal of
this section is to obtain such an estimate.  We proceed in two steps.
First, we use native filtering arguments as in \cite{BaxChiLip} to obtain
an a.s.\ exponential estimate for $|D\pi_{0,t}(\nu)|$.  As the laws of the
observation processes generated by signals with different initial
distributions and jump rates are equivalent, we can extend this a.s.\
bound to $|D\tilde\pi_{s,t}(\mu)|$.  We find, however, that the proportionality
constant in the exponential estimate depends on $\mu$ and diverges as
$\mu$ approaches the boundary of the simplex.  This makes a pathwise bound
on $|D\tilde\pi_{s,t}(\pi_s)|$ difficult to obtain, as $\pi_s$ can get
arbitrarily close to the boundary of the simplex on the infinite time
interval.  Instead, we proceed to find a uniform bound on
$\mathbf{E}|D\tilde\pi_{s,t}(\pi_s)|$.

We begin by recalling a few useful results from \cite{BaxChiLip}.

\begin{lem}
Assume $\mu,\nu$ are in the interior of the simplex.  Then
\begin{equation}
\label{eq:wronginirep}
\pi_t^i(\mu)=\frac{\sum_j (\mu^j/\nu^j)\mathbf{P}(X_0=a_j,
X_t=a_i|\mathscr{F}_t^Y)}{\sum_j(\mu^j/\nu^j)\mathbf{P}(X_0=a_j|\mathscr{F}_t^Y)}.
\end{equation}
\end{lem}

\begin{pf}
Define a new measure $\mathbf{P}^\mu\sim\mathbf{P}$ through
\[
\frac{d\mathbf{P}^\mu}{d\mathbf{P}}=\frac{d\mu}{d\nu}(X_0).
\]
It is not difficult to verify that under $\mathbf{P}^\mu$, $X_t$ is still
a finite-state Markov process with intensities matrix $\Lambda$ but with
initial distribution $\mathbf{P}^\mu(X_0=a_i)=\mu^i$.  Hence evidently
$\pi_t^i(\mu)=\mathbf{P}^\mu(X_t=a_i|\mathscr{F}_t^Y)$.  Using the
usual change of measure formula for conditional expectations, we
can write
\[
\pi_t^i(\mu)=\mathbf{E}_{\mathbf{P}^\mu}(I_{X_t=a_i}|\mathscr{F}_t^Y)=
\frac{\mathbf{E}(I_{X_t=a_i}\,(d\mu/d\nu)(X_0)|\mathscr{F}_t^Y)}
{\mathbf{E}((d\mu/d\nu)(X_0)|\mathscr{F}_t^Y)}.
\]
The result now follows immediately.
\end{pf}

For the proof of the following lemma we refer to
\cite{BaxChiLip}, Lemma 5.7, page 662.

\begin{lem}
\label{lem:smoothingestimate}
Define $\rho^{ji}_t=\mathbf{P}(X_0=a_j|\mathscr{F}_t^Y,\,X_t=a_i)$.
Assume that $\lambda_{ij}>0$ $\forall i\ne j$.
Then for any $t\ge 0$ we have the a.s. bound
\[
\max_{j,k,\ell}|\rho^{jk}_t-\rho^{j\ell}_t|\le\exp\biggl(
-2t\min_{p,q\ne p}\sqrt{\lambda_{pq}\lambda_{qp}}\biggr).
\]
\end{lem}

We are now ready to obtain some useful estimates.

\begin{prop}
\label{pro:jacoboundsimple}
Let $\lambda_{ij}>0$ $\forall i\ne j$ and $\nu\in\mathcal{S}^{d-1}$, $v\in
T\mathcal{S}^{d-1}$.  Then a.s.
\[
|D\pi_t(\nu)\cdot v|\le\sum_k\frac{|v^k|}{\nu^k}\exp\biggl(
-2t\min_{p,q\ne p}\sqrt{\lambda_{pq}\lambda_{qp}}\biggr).
\]
\end{prop}

\begin{pf}
We can calculate directly the directional derivative of (\ref{eq:wronginirep}):
\[
\bigl(D\pi_t(\mu)\cdot v\bigr)^i= \frac{\sum_j {(v^j/\nu^j)}(\mathbf{P}(X_0=a_j,
X_t=a_i|\mathscr{F}_t^Y)-\pi_t^i(\mu)\mathbf{P}(X_0=a_j|\mathscr{F}_t^Y))}
{\sum_j{(\mu^j/\nu^j)}\,\mathbf{P}(X_0=a_j|\mathscr{F}_t^Y)}.
\]
Setting $\mu=\nu$, we obtain after some simple manipulations
\[
\bigl(D\pi_t(\nu)\cdot v\bigr)^i=
\pi_t^i(\nu)\sum_{j,k} (v^j/\nu^j)\,\pi_t^k(\nu)\,(\rho_t^{ji}-\rho_t^{jk}).
\]
The result follows from Lemma \ref{lem:smoothingestimate}.
\end{pf}

To obtain this bound we had to use the true initial distribution $\nu$,
jump rates $\lambda_{ij}$ and observation function $h$.  However, the
almost sure nature of the result allows us to drop these requirements.

\begin{cor}
\label{cor:forgetting}
Let $\tilde\lambda_{ij}>0$ $\forall i\ne j$ and $\mu\in\mathcal{S}^{d-1}$,
$v\in T\mathcal{S}^{d-1}$.  Then a.s.
\begin{equation}
\label{eq:forgetting}
|D\tilde\pi_{s,t}(\mu)\cdot v|\le \sum_k\frac{|v^k|}{\mu^k}
\exp\biggl(-2(t-s)\min_{p,q\ne p}\sqrt{\tilde\lambda_{pq}
\tilde\lambda_{qp}}\biggr).
\end{equation}
Moreover, the result still holds if $\mu,v$ are
$\mathscr{F}_s^Y$-measurable random variables with values a.s. in
$\mathcal{S}^{d-1}$ and $T\mathcal{S}^{d-1}$, respectively.
\end{cor}

\begin{pf}
Note that we can write
$\tilde\pi_{0,t}^i(\mu)=\mathbf{\tilde P}^\mu(X_t=a_i|\mathscr{F}_t^Y)$,
where $\mathbf{\tilde P}^\mu$ is the
measure under which $X_t$ has transition intensities matrix $\tilde\Lambda$ and
initial distribution~$\mu$, and $dY_t = \mathrm{\tilde h}(X_t)\,dt+d\tilde B_t$
where $\tilde B_t$ is a Wiener process independent of~$X_t$. But
$\mathbf{\tilde P}^\mu$ and $\mathbf{P}$ are equivalent measures (by the
Girsanov theorem and \cite{rogersw}, Section~IV.22), so that the result
for $s=0$ follows trivially from Proposition \ref{pro:jacoboundsimple}.
The result for $s>0$ follows directly as the Wonham equation is time homogeneous.

To show that the result still holds when $\mu,v$ are random, note that
$\tilde\pi_{s,t}$ only depends on the observation increments in the interval
$[s,t]$, that is, $D\tilde\pi_{s,t}(\mu)\cdot v$ is
$\mathscr{F}_{[s,t]}^Y$-measurable where $\mathscr{F}_{[s,t]}^Y=
\sigma\{Y_r-Y_s\dvtx s\le r\le t\}$.  Under the equivalent measure $\mathbf{Q}$
introduced in Section \ref{sec:prelims}, $Y$ is a Wiener process and hence
$\mathscr{F}_{[s,t]}^Y$ and $\mathscr{F}_{s}^Y$ are independent.  It
follows from the bound with constant $\mu,v$ that
\[
\mathbf{E_Q}\bigl(I_{|D\tilde\pi_{s,t}(\mu)\cdot v|\le(*)}|\sigma\{\mu,v\}\bigr)=1,
\qquad \mathbf{Q}\mbox{-a.s.},
\]
where $(*)$ is the right-hand side of (\ref{eq:forgetting}).
Hence $\mathbf{E_Q}(I_{|D\tilde\pi_{s,t}(\mu)\cdot v|\le(*)})=1$, and the
statement follows from $\mathbf{P}\sim\mathbf{Q}$.
\end{pf}

\begin{prop}
\label{pro:forgetting}
Let $\tilde\lambda_{ij}>0$ $\forall i\ne j$ and
$\mu_1,\mu_2\in\mathcal{S}^{d-1}$.  Then a.s.
\[
|\tilde\pi_{s,t}(\mu_2)-\tilde\pi_{s,t}(\mu_1)|\le C|\mu_2-\mu_1|
\exp\biggl(-2(t-s)\min_{p,q\ne p}\sqrt{\tilde\lambda_{pq}
\tilde\lambda_{qp}}\biggr),
\]
where $C=\max\{1/\mu_1^k,1/\mu_2^k\dvtx k=1,\ldots,d\}$.
\end{prop}

\begin{pf}
Define $\gamma(u)=\tilde\pi_{s,t}(\mu_1+u(\mu_2-\mu_1))$, $u\in[0,1]$. Then
\[
\tilde\pi_{s,t}(\mu_2)-\tilde\pi_{s,t}(\mu_1)=\int_0^1 \frac{d\gamma}{du}\,du=
\int_0^1 D\tilde\pi_{s,t}\bigl(\mu_1+u(\mu_2-\mu_1)\bigr)\cdot(\mu_2-\mu_1)\,du.
\]
We can thus estimate
\[
|\tilde\pi_{s,t}(\mu_2)-\tilde\pi_{s,t}(\mu_1)|\le
\sup_{u\in[0,1]}\big|D\tilde\pi_{s,t}\bigl(\mu_1+u(\mu_2-\mu_1)\bigr)\cdot
(\mu_2-\mu_1)\big|.
\]
The result now follows from Corollary \ref{cor:forgetting}.
\end{pf}

Corollary \ref{cor:forgetting} and Proposition \ref{pro:forgetting} are
exactly what we need to establish boundedness of equation
(\ref{eq:simplebound}).  Note, however, that the right-hand side of
(\ref{eq:forgetting}) is proportional to $1/\mu^i$, and we must estimate
$|D\tilde\pi_{s,t}(\pi_s)|$. Though we established in Section \ref{sec:prelims}
that $\pi_s$ cannot hit the boundary of the simplex in finite time, it can
get arbitrarily close to the boundary during the infinite time interval,
thus rendering the right-hand side of equation (\ref{eq:forgetting})
arbitrarily large.  If we can establish that $\sup_{s\ge
0}\mathbf{E}(1/\min_k\pi_s^k)<\infty$, however, then we can control
$\mathbf{E}|D\tilde\pi_{s,t}(\pi_s)|$ to obtain a useful bound.

We begin with an auxiliary integrability property of $\pi_t$:

\begin{lem}
\label{lem:finiteexpect}
Let $\nu\in\mathcal{S}^{d-1}$ and $T<\infty$. Then
\[
{\bf E} \int_0^T (\pi_s^i)^{-k}\,ds<\infty\qquad
\forall\,i=1,\ldots,d, k\ge 1.
\]
\end{lem}

\begin{pf}
Applying It\^o's rule to the Wonham equation gives
\[
d\log\pi_t^i=\biggl(\lambda_{ii}-\frac{1}{2}(h^i-h^*\pi_t)^2\biggr)\,dt
+\sum_{j\ne i}\lambda_{ji}\frac{\pi_t^j}
{\pi_t^i}\,dt+(h^i-h^*\pi_t)\,dW_t,
\]
where the innovation $dW_t=dY_t-h^*\pi_t\,dt$ is an
$\mathscr{F}_t^Y$-Wiener process.  The application of It\^o's rule is
justified by a standard localization argument, as $\pi_t$ is in
$\mathcal{S}^{d-1}$ for all $t\ge 0$ a.s.\ and $\log x$ is smooth in
$(0,1)$.  As $\lambda_{ij}\ge 0$ for $j\ne i$, we estimate
\[
-k\log\pi_t^i\le -k\log\nu^i
-k\lambda_{ii}t + \frac{k}{2}\max_j(h^i-h^j)^2\,t
-k\int_0^t (h^i-h^*\pi_s)\,dW_s.
\]
But as $h^i-h^*\pi_t$ is bounded, Novikov's condition is satisfied and hence
\[
\mathbf{E}\exp\biggl(-k\int_0^t (h^i-h^*\pi_s)\,dW_s
-\frac{k^2}{2}\int_0^t (h^i-h^*\pi_s)^2\,ds\biggr)=1.
\]
Estimating the time integral, we obtain
\[
\mathbf{E}(\pi_t^i)^{-k} \le (\nu^i)^{-k}
\exp\biggl(-k\lambda_{ii}t + \tfrac{1}{2}k(k+1)\max_j(h^i-h^j)^2\,t\biggr).
\]
The lemma now follows by the Fubini--Tonelli theorem, as
$(\pi_s^i)^{-k}\ge 0$ a.s.
\end{pf}

We are now in a position to bound $\sup_{t\ge 0}{\bf E}(1/\min_i\pi_t^i)$.

\begin{prop}
\label{prop:inf}
Let $\nu\in\mathcal{S}^{d-1}$ and suppose that $\lambda_{ij}>0$ $\forall i\ne j$. Then
\[
\sup_{t\ge 0}{\bf E}\biggl(\frac{1}{\min_i\pi_t^i}\biggr)<\infty.
\]
\end{prop}

\begin{pf}
By It\^o's rule and using the standard localization argument, we obtain
\begin{eqnarray*}
(\pi_t^i)^{-1} &=& (\nu^i)^{-1}
-\int_0^t\lambda_{ii}(\pi_s^i)^{-1}\,ds
-\int_0^t(\pi_s^i)^{-2}\sum_{j\ne i}\lambda_{ji}\pi_s^j\,ds
\\
&&{}-\int_0^t(\pi_s^i)^{-1}(h^i-h^*\pi_s)\,dW_s
+\int_0^t(\pi_s^i)^{-1}(h^i-h^*\pi_s)^2\,ds,
\end{eqnarray*}
where $W_t$ is the innovations Wiener process.  Using Lemma
\ref{lem:finiteexpect} we find
\[
\mathbf{E}\int_0^t (\pi_s^i)^{-2}(h^i-h^*\pi_s)^2\,ds\le
\max_j(h^i-h^j)^2\mathbf{E}\int_0^t (\pi_s^i)^{-2}\,ds<\infty,
\]
so the expectation of the stochastic integral term vanishes.  Using the
Fubini--Tonelli theorem, we can thus write
\begin{eqnarray*}
\mathbf{E}((\pi_t^i)^{-1}) &=& (\nu^i)^{-1} -\int_0^t\lambda_{ii}
\mathbf{E}((\pi_s^i)^{-1})\,ds
\\
&&{}-\int_0^t\mathbf{E}\Biggl((\pi_s^i)^{-2}\sum_{j\ne i}\lambda_{ji}\pi_s^j\Biggr)\,ds
+\int_0^t\mathbf{E}\bigl((\pi_s^i)^{-1}(h^i-h^*\pi_s)^2\bigr)\,ds.
\end{eqnarray*}
Taking the derivative and estimating each of the terms, we obtain
\[
\frac{dM_t^i}{dt}\le-\min_{j\ne i}\lambda_{ji}\,(M_t^i)^2
+\biggl(|\lambda_{ii}|+\min_{j\ne i}\lambda_{ji}+\max_j(h^i-h^j)^2 \biggr) M_t^i,
\]
where we have written $M_t^i = \mathbf{E}((\pi_t^i)^{-1})$ and we have
used $(M_t^i)^2\le \mathbf{E}(\pi_t^i)^{-2}$ by Jensen's inequality.
Using the estimate
\[
-K_1^i(M_t^i)^2+K_2^iM_t^i\le -K_2^iM_t^i+\frac{(K_2^i)^2}{K_1^i}
\qquad  \mbox{for } K_1^i>0,
\]
we now obtain
\[
\frac{dM_t^i}{dt}\le K_2^i\biggl(\frac{K_2^i}{K_1^i}-M_t^i\biggr),
\qquad  K_2^i=|\lambda_{ii}|+\min_{j\ne i}\lambda_{ji}+\max_j(h^i-h^j)^2,
\]
where $K_1^i=\min_{j\ne i}\lambda_{ji}>0$.  Consequently we obtain
\[
M_t^i\le e^{-K_2^i t}(\nu^i)^{-1}+\frac{(K_2^i)^2}{K_1^i}
e^{-K_2^i t}\int_0^t e^{K_2^i s}\,ds=
e^{-K_2^i t}(\nu^i)^{-1}+\frac{K_2^i}{K_1^i}(1-e^{-K_2^i t}).
\]
We can now estimate
\[
\sup_{t\ge 0}{\bf E}\biggl(\frac{1}{\min_i\pi_t^i}\biggr)
\le \sum_{i=1}^d\sup_{t\ge 0}{\bf E}\biggl(\frac{1}{\pi_t^i}\biggr)
\le \sum_{i=1}^d\biggl(\frac{1}{\nu^i}\vee \frac{K_2^i}{K_1^i}\biggr)<\infty,
\]
which is what we set out to prove.
\end{pf}

We can now prove Theorem \ref{thm:mainresult} for the special case
$\tilde h=h$. Using equation (\ref{eq:simplebound}), Corollary
\ref{cor:forgetting}, Proposition \ref{pro:forgetting} and Proposition
\ref{prop:inf}, we obtain
\begin{eqnarray*}
&& \mathbf{E}|\pi_t-\tilde\pi_t(\mu)|
\\
&&\qquad \le |\mu-\nu|\max_k\biggl\{\frac{1}{\mu^k}\vee
\frac{1}{\nu^k}\biggr\}
\exp\biggl(-2t\min_{p,q\ne p}\sqrt{\tilde\lambda_{pq}
\tilde\lambda_{qp}}\biggr)
\\
&&\qquad\phantom{\leq}{}+|\Lambda^*-\tilde\Lambda^*|
\sup_{s\ge 0}\mathbf{E}\biggl(1/\min_k\pi_s^k\biggr)
\int_0^t\exp\biggl(-2(t-s)\min_{p,q\ne p}\sqrt{\tilde\lambda_{pq}
\tilde\lambda_{qp}}\biggr)\,ds,
\end{eqnarray*}
where $|\Lambda^*-\tilde\Lambda^*|=\sup\{|(\Lambda^*-\tilde\Lambda^*)\mu|\dvtx
\mu\in\mathcal{S}^{d-1}\}$. Thus
\[
\mathbf{E}|\pi_t-\tilde\pi_t(\mu)|\le
|\mu-\nu|\max_k\biggl\{\frac{1}{\mu^k}\vee
\frac{1}{\nu^k}\biggr\}e^{-\beta t}
+|\Lambda^*-\tilde\Lambda^*|
\frac{\sup_{s\ge 0}\mathbf{E}(1/\min_k\pi_s^k)}{\beta},
\]
where we have written
$\beta=2\min_{p,q\ne p}(\tilde\lambda_{pq}\tilde\lambda_{qp})^{1/2}$.
The result follows directly using $\|\pi_t-\tilde\pi_t(\mu)\|^2\le
|\pi_t-\tilde\pi_t(\mu)|$ [as $|\pi_t^i-\tilde\pi_t(\mu)^i|\le 1$].

\section{Model robustness of the Wonham filter}
\label{sec:robustness}

We are now ready to proceed to the general case where the initial density,
the transition intensities matrix and the observation function can all be
misspecified.  The simplicity of the special case $\tilde h=h$ that we have
treated up to this point is due to the fact that in the calculation of
equation (\ref{eq:itodropsout}), the stochastic integral term drops out and we
can proceed with the calculation using only ordinary calculus.  In the
general case we cannot get rid of the stochastic integral, and hence we
run into anticipativity problems in the next step of the calculation.

We solve this problem by using anticipative stochastic integrals in the
sense of Skorokhod, rather than the usual It\^o integral (which is a
special case of the Skorokhod integral defined for adapted processes
only).  Though the Skorokhod integral is more general than the It\^o
integral in the sense that it allows some anticipating integrands, it is
less general in that we have to integrate against a Wiener process (rather
than against an arbitrary semimartingale), and that the integrands should
be functionals of the driving Wiener process.  In our setup, the most
convenient way to deal with this is to operate exclusively under the
measure $\mathbf{Q}$ of Section~\ref{sec:prelims}, under which the
observation process $Y$ is a Wiener process.  At the end of the day we can
calculate the relevant expectation with respect to the measure $\mathbf{P}$
by using the explicit expression for the Radon--Nikodym derivative
$d\mathbf{P}/d\mathbf{Q}$.  The fact that the integrands must be
functionals of the underlying Wiener process is not an issue, as both the
approximate and exact filters are functionals of the observations only.

Our setup is further detailed in Appendix \ref{sec:malliavin}, together
with a review of the relevant results from the Malliavin calculus and
anticipative stochastic calculus.  Below we will use the notation and
results from this appendix without further comment.  We will also refer to
Appendix \ref{sec:proofs} for some results on smoothness of the various
integrands we encounter; these results are not central to the
calculations, but are required for the application of the theory in
Appendix \ref{sec:malliavin}.

We begin by obtaining an anticipative version of Proposition
\ref{pro:simplebound}.  Note that this result is precisely of the form one
would expect.  The first two lines follow the formula for the distance
between two flows as one would guess, for example, from the discussion in the
\hyperref[sec:intro]{Introduction}; the last line is an It\^o correction term which contains
second derivatives of the filter with respect to its initial condition.

\begin{prop}
\label{pro:anticipativediff}
The difference between $\pi_t$ and $\tilde\pi_t$ satisfies
\begin{eqnarray*}
\pi_t-\tilde\pi_t &=&\int_0^t D\tilde\pi_{r,t}(\pi_r)\cdot\Delta_\Lambda\pi_r\,dr
+\int_0^t D\tilde\pi_{r,t}(\pi_r)\cdot\Delta_H(\pi_r)\,dY_r
\\
&&{} -\int_0^t D\tilde\pi_{r,t}(\pi_r)\cdot[h^*\pi_r\,(H-h^*\pi_r)\pi_r-
\tilde h^*\pi_r\,(\tilde H-\tilde h^*\pi_r)\pi_r]\,dr
\\
&&{} +\tfrac{1}{2}\int_0^t[D^2\tilde\pi_{r,t}(\pi_r)\cdot(H-h^*\pi_r)\pi_r
-D^2\tilde\pi_{r,t}(\pi_r)\cdot(\tilde H-\tilde h^*\pi_r)\pi_r]\,dr,
\end{eqnarray*}
where the stochastic integral is a Skorokhod integral and we have written
$\Delta_\Lambda=\Lambda^*-\tilde\Lambda^*$,
$\Delta_H(\pi)=(H-h^*\pi)\pi-(\tilde H-\tilde h^*\pi)\pi$, and
$D^2\tilde\pi_{r,t}(\mu)\cdot v$ is the directional derivative of
$D\tilde\pi_{r,t}(\mu)\cdot v$ with respect to $\mu\in\mathcal{S}^{d-1}$ in the
direction $v\in T\mathcal{S}^{d-1}$.
\end{prop}

\begin{pf}
Fix some $T>t$.  We begin by evaluating, using It\^o's rule and equation
(\ref{eq:zakaiinv}),
\begin{eqnarray*}
\tilde U_{0,s}^{-1}U_{0,s}\nu &=& \nu +\int_0^s\tilde U_{0,r}^{-1}
(\Lambda^*-\tilde\Lambda^*)U_{0,r}\nu\,dr
\\
&&{} -\int_0^s\tilde U_{0,r}^{-1}\tilde H(H-\tilde H)U_{0,r}\nu\,dr
+\int_0^s\tilde U_{0,r}^{-1}(H-\tilde H)U_{0,r}\nu\,dY_r.
\end{eqnarray*}
Now multiply from the left by $\tilde U_{0,t}$; we wish to use Lemma
\ref{lem:bringintoskor} to bring $\tilde U_{0,t}$ into the Skorokhod integral
term, that is, we claim that
\begin{eqnarray*}
\tilde U_{s,t}U_{0,s}\nu &=& \tilde U_{0,t}\nu
+\int_0^s\tilde U_{r,t}(\Lambda^*-\tilde\Lambda^*)U_{0,r}\nu\,dr
-\int_0^s\tilde U_{r,t}\tilde H(H-\tilde H)U_{0,r}\nu\,dr
\\
&&{} +\int_0^s\tilde U_{r,t}(H-\tilde H)U_{0,r}\nu\,dY_r
+\int_0^s(\mall_r\tilde U_{0,t})\tilde U_{0,r}^{-1}(H-\tilde H)U_{0,r}\nu\,dr.
\end{eqnarray*}
To justify this expression we need to verify the integrability conditions
of Lemma~\ref{lem:bringintoskor}.  Note that all matrix
elements of $\tilde U_{s,t}$ are in $\mathbb{D}^\infty$ $\forall 0\le s\le
t<T$, and that
\[
\mall_r\tilde U_{s,t}=\cases{
0, & \quad\mbox{a.e. }$r\notin[s,t]$, \cr
\tilde U_{r,t}\tilde H\tilde U_{s,r}, & \quad\mbox{a.e. }$r\in[s,t]$.}
\]
This follows directly from Proposition \ref{pro:sdemalliavin} and Lemma
\ref{lem:zakaiflow} (note that the same result holds for $U_{s,t}$ if we
replace $\tilde H$ by $H$ and $\tilde U$ by $U$).  Once we plug this result into
the expression above, the corresponding integrability conditions can be
verified explicitly, see Lemma \ref{lem:msqintg}, and hence we have
verified that
\[
\tilde U_{s,t}U_{0,s}\nu=\tilde U_{0,t}\nu
+\int_0^s\tilde U_{r,t}(\Lambda^*-\tilde\Lambda^*)U_{0,r}\nu\,dr
+\int_0^s\tilde U_{r,t}(H-\tilde H)U_{0,r}\nu\,dY_r.
\]
Next we would like to apply the anticipating It\^o rule, Proposition
\ref{pro:anticipatingito}, with the function
$\Sigma\dvtx\mathbb{R}^d_{++}\to\mathcal{S}^{d-1}$, $\Sigma(x)=x/|x|$.
To this end we have to verify a set of technical conditions, see Lemma
\ref{lem:canapplyito}.  We obtain
\begin{eqnarray*}
&& \Sigma(\tilde U_{s,t}U_{0,s}\nu)
\\
&&\qquad =\Sigma(\tilde U_{0,t}\nu)
+\int_0^s D\Sigma(\tilde U_{r,t}U_{0,r}\nu)
\tilde U_{r,t}(\Lambda^*-\tilde\Lambda^*)U_{0,r}\nu\,dr
\\
&&\qquad\phantom{=}{} +\frac{1}{2}\sum_{k,\ell}\int_0^s\frac{\partial^2\Sigma}
{\partial x^k\,\partial x^\ell}
(\tilde U_{r,t}U_{0,r}\nu)(\nab_r\tilde U_{r,t}U_{0,r}\nu)^k
\bigl(\tilde U_{r,t}(H-\tilde H)U_{0,r}\nu\bigr)^\ell\,dr
\\
&&\qquad\phantom{=}{} +\int_0^s D\Sigma(\tilde U_{r,t}U_{0,r}\nu) \tilde U_{r,t}(H-\tilde H)U_{0,r}\nu\,dY_r.
\end{eqnarray*}
We need to evaluate $\nab_r\tilde U_{r,t}U_{0,r}\nu$.  Using
Proposition \ref{pro:chainrulesmooth}, we calculate
\[
\lim_{\varepsilon\searrow 0}\mall_r\tilde U_{r+\varepsilon,t}U_{0,r+\varepsilon}\nu =
\lim_{\varepsilon\searrow 0}\tilde U_{r+\varepsilon,t}U_{r,r+\varepsilon}HU_{0,r}\nu =
\tilde U_{r,t}HU_{0,r}\nu,
\]
and similarly
\[
\lim_{\varepsilon\searrow 0}\mall_r\tilde U_{r-\varepsilon,t}U_{0,r-\varepsilon}\nu =
\lim_{\varepsilon\searrow 0}\tilde U_{r,t}\tilde H\tilde
U_{r-\varepsilon,r}U_{0,r-\varepsilon}\nu =\tilde U_{r,t}\tilde H U_{0,r}\nu.
\]
After some rearranging, we obtain
\begin{eqnarray*}
\Sigma(\tilde U_{s,t}U_{0,s}\nu) &=& \Sigma(\tilde U_{0,t}\nu)
+\int_0^s D\Sigma(\tilde U_{r,t}U_{0,r}\nu)
\tilde U_{r,t}(\Lambda^*-\tilde\Lambda^*)U_{0,r}\nu\,dr
\\
&&{} +\frac{1}{2}\sum_{k,\ell}\int_0^s\frac{\partial^2\Sigma}{\partial x^k\,\partial x^\ell}
(\tilde U_{r,t}U_{0,r}\nu)(\tilde U_{r,t}HU_{0,r}\nu)^k
(\tilde U_{r,t}HU_{0,r}\nu)^\ell\,dr
\\
&&{}-\frac{1}{2}\sum_{k,\ell} \int_0^s\frac{\partial^2\Sigma}{\partial x^k\,\partial x^\ell}
(\tilde U_{r,t}U_{0,r}\nu)(\tilde U_{r,t}\tilde H U_{0,r}\nu)^k
(\tilde U_{r,t}\tilde H U_{0,r}\nu)^\ell\,dr
\\
&&{} +\int_0^s D\Sigma(\tilde U_{r,t}U_{0,r}\nu)
\tilde U_{r,t}(H-\tilde H)U_{0,r}\nu\,dY_r.
\end{eqnarray*}
From this point onward we will set $s=t$.  We will need (on $\mathbb{R}^d_{++}$)
\[
D^2\Sigma^{ik\ell}(x)=\frac{\partial^2\Sigma^i(x)}{\partial x^k\,\partial x^\ell}=
 -\frac{1}{|x|}\bigl(D\Sigma^{ik}(x)+D\Sigma^{i\ell}(x)\bigr).
\]
Recall that $D\Sigma(\alpha x)=\alpha^{-1}D\Sigma(x)$; it follows that also
$D^2\Sigma(\alpha x)=\alpha^{-2}D^2\Sigma(x)$ for $\alpha>0$.  Using these
expressions with $\alpha=|U_{0,r}\nu|$, we get
\begin{eqnarray*}
\pi_t-\tilde\pi_t &=& \int_0^t D\Sigma(\tilde U_{r,t}\pi_r)\tilde U_{r,t}\Delta_\Lambda\pi_r\,dr
+\int_0^t D\Sigma(\tilde U_{r,t}\pi_r) \tilde U_{r,t}(H-\tilde H)\pi_r\,dY_r
\\
&&{} +\frac{1}{2}\sum_{k,\ell}\int_0^t\frac{\partial^2\Sigma}{\partial x^k\,\partial x^\ell}
(\tilde U_{r,t}\pi_r)(\tilde U_{r,t}H\pi_r)^k(\tilde U_{r,t}H\pi_r)^\ell\,dr
\\
&&{}-\frac{1}{2}\sum_{k,\ell} \int_0^t\frac{\partial^2\Sigma}{\partial x^k\,\partial x^\ell}
(\tilde U_{r,t}\pi_r)(\tilde U_{r,t}\tilde H\pi_r)^k(\tilde U_{r,t}\tilde H\pi_r)^\ell\,dr.
\end{eqnarray*}
Next we want to express the integrands in terms of $D\tilde\pi_{r,t}(\pi_r)\cdot v$,
and so on, rather than in terms of $D\Sigma(x)$.  Recall that
$D\tilde\pi_{r,t}(\pi_r)\cdot v=D\Sigma(\tilde U_{r,t}\pi_r)\tilde U_{r,t}v$ when
$v\in T\mathcal{S}^{d-1}$.  Similar terms appear in the expression above,
but, for example, $\tilde H\pi_r\notin T\mathcal{S}^{d-1}$.  To rewrite the
expression in the desired form, we use that
$D\Sigma(\tilde U_{r,t}\pi_r)\tilde U_{r,t}\pi_r=0$.  Hence
\begin{eqnarray*}
D\Sigma(\tilde U_{r,t}\pi_r)\tilde U_{r,t}\tilde H\pi_r &=&
D\Sigma(\tilde U_{r,t}\pi_r)\tilde U_{r,t}(\tilde H-\tilde h^*\pi_r)\pi_r
\\
&=& D\tilde\pi_{r,t}(\pi_r)\cdot (\tilde H-\tilde h^*\pi_r)\pi_r
\end{eqnarray*}
and similarly for the other terms.  Note also that
\[
\sum_k D^2\Sigma^{ik\ell}(\tilde U_{r,t}\pi_r)(\tilde U_{r,t}\pi_r)^k=
 -D\Sigma^{i\ell}(\tilde U_{r,t}\pi_r).
\]
Substituting this into the expression for $\pi_t-\tilde\pi_t$ and rearranging,
we obtain
\begin{eqnarray*}
&& \pi_t-\tilde\pi_t
\\
&&\qquad = \int_0^t D\tilde\pi_{r,t}(\pi_r)\cdot\Delta_\Lambda\pi_r\,dr
+\int_0^t D\tilde\pi_{r,t}(\pi_r)\cdot\Delta_H(\pi_r)\,dY_r
\\
&&\hspace*{-1pt}\qquad\phantom{=}{} -\int_0^t D\tilde\pi_{r,t}(\pi_r)\cdot[h^*\pi_r\,(H-h^*\pi_r)\pi_r-
\tilde h^*\pi_r\,(\tilde H-\tilde h^*\pi_r)\pi_r]\,dr
\\
&&\hspace*{-1pt}\qquad\phantom{=}{} +\frac{1}{2}\sum_{k,\ell}\int_0^t\frac{\partial^2\Sigma}
{\partial x^k\,\partial x^\ell}
(\tilde U_{r,t}\pi_r)\bigl(\tilde U_{r,t}(H-h^*\pi_r)\pi_r\bigr)^k
\bigl(\tilde U_{r,t}(H-h^*\pi_r)\pi_r\bigr)^\ell\,dr
\\
&&\hspace*{-1pt}\qquad\phantom{=}{} -\frac{1}{2}\sum_{k,\ell}\int_0^t\frac{\partial^2\Sigma}
{\partial x^k\,\partial x^\ell}
(\tilde U_{r,t}\pi_r)\bigl(\tilde U_{r,t}(\tilde H-\tilde h^*\pi_r)\pi_r\bigr)^k
\bigl(\tilde U_{r,t}(\tilde H-\tilde h^*\pi_r)\pi_r\bigr)^\ell\,dr.
\end{eqnarray*}
It remains to note that we can write
\[
\bigl(D^2\tilde\pi_{s,t}(\mu)\cdot v\bigr)^i =
\sum_{k,\ell}D^2\Sigma^{ik\ell}(\tilde U_{s,t}\mu)
(\tilde U_{s,t}v)^k(\tilde U_{s,t}v)^\ell.
\]
The result follows immediately.
\end{pf}

\begin{rem}
We have allowed misspecification of most model parameters of the Wonham
filter.  One exception is the observation noise intensity: we have not
considered observations of the form $dY_t=\mathrm{h}(X_t)\,dt+\sigma\,dB_t$
with $\sigma\ne 1$; in other words, the quadratic variation of $Y_t$ is
assumed to be known $[Y,Y]_t = t$.  We do not consider this a significant
drawback as the quadratic variation can be determined directly from the
observation process $Y_t$. On the other hand, the model parameters
$\nu,\Lambda,h$ are ``hidden'' and would have to be estimated, making
these quantities much more prone to modeling errors.

If we allow misspecification of $\sigma$, we would have to be careful to
specify in which way the filter is implemented: in this case, the
normalized solution of the misspecified Zakai equation no longer coincides
with the solution of the misspecified Wonham equation.  Hence one obtains
a different error estimate depending on whether the normalized
solution of the misspecified Zakai equation, or the solution of the
misspecified Wonham equation, is compared to the exact filter.  Both
cases can be treated using similar methods, but we do not pursue this here.
\end{rem}

Let $e_t=\pi_t-\tilde\pi_t$.  We wish to estimate the norm of
$e_t$.  Unfortunately, we can no longer use the triangle inequality as in
Section \ref{sec:prelims} due to the presence of the stochastic integral;
instead, we choose to calculate $\|e_t\|^2$, which is readily estimated.

\begin{lem}
\label{lem:filtererror}
The filtering error can be estimated by
\begin{eqnarray*}
&&\mathbf{E_P}\|e_t\|^2
\\
&&\qquad \le
\int_0^t \mathbf{E_P}|D\tilde\pi_{r,t}(\pi_r)\cdot\Delta_\Lambda\pi_r|\,dr
+ K\int_0^t\mathbf{E_P}|D\tilde\pi_{r,t}(\pi_r)\cdot\Delta_H(\pi_r)|\,dr
\\
&&\qquad\phantom{\le}{}+\int_0^t \mathbf{E_P}\big|D\tilde\pi_{r,t}(\pi_r)\cdot
\bigl(h^*\pi_r\,(H-h^*\pi_r)\pi_r-
\tilde h^*\pi_r\,(\tilde H-\tilde h^*\pi_r)\pi_r\bigr)\big|\,dr
\\
&&\qquad\phantom{\le}{} +\tfrac{1}{2}\int_0^t \mathbf{E_P}|
D^2\tilde\pi_{r,t}(\pi_r)\cdot(H-h^*\pi_r)\pi_r
-D^2\tilde\pi_{r,t}(\pi_r)\cdot(\tilde H-\tilde h^*\pi_r)\pi_r|\,dr,
\end{eqnarray*}
where $K=2\max_k|h^k|+\max_k|\tilde h^k|$.
\end{lem}

\begin{pf}
We wish to calculate $\mathbf{E_P}\|e_t\|^2=\mathbf{E_P}e_t^*e_t$.  Using
Proposition \ref{pro:anticipativediff}, we obtain
\begin{eqnarray*}
&&\mathbf{E_P}\|e_t\|^2
\\
&&\qquad = \int_0^t \mathbf{E_P}
e_t^*D\tilde\pi_{r,t}(\pi_r)\cdot\Delta_\Lambda\pi_r\,dr
\\
&&\qquad\phantom{=}{}+ \mathbf{E_P}\biggl[e_t^*
\int_0^t D\tilde\pi_{r,t}(\pi_r)\cdot\Delta_H(\pi_r)\,dY_r\biggr]
\\
&&\qquad\phantom{=}{} -\int_0^t \mathbf{E_P}\,e_t^*D\tilde\pi_{r,t}(\pi_r)\cdot
[h^*\pi_r\,(H-h^*\pi_r)\pi_r- \tilde h^*\pi_r\,(\tilde H-\tilde h^*\pi_r)\pi_r]\,dr
\\
&&\qquad\phantom{=}{} +\tfrac{1}{2}\int_0^t \mathbf{E_P}
e_t^*[D^2\tilde\pi_{r,t}(\pi_r)\cdot(H-h^*\pi_r)\pi_r
\\
&&\hspace*{90pt}{}-D^2\tilde\pi_{r,t}(\pi_r)\cdot(\tilde H-\tilde h^*\pi_r)\pi_r]\,dr.
\end{eqnarray*}
The chief difficulty is the stochastic integral term.
Using equation (\ref{eq:refmeasure}), we can write
\begin{eqnarray*}
&& \mathbf{E_P}\biggl[e_t^* \int_0^t D\tilde\pi_{r,t}(\pi_r)\cdot\Delta_H(\pi_r)\,dY_r\biggr]
\\
&&\qquad = \mathbf{E_Q}\biggl[|U_{0,t}\nu|\,e_t^*
\int_0^t D\tilde\pi_{r,t}(\pi_r)\cdot\Delta_H(\pi_r)\,dY_r\biggr].
\end{eqnarray*}
We would like to apply equation (\ref{eq:malladjoint}) to evaluate this
expression. First, we must establish that the integrand is in
$\operatorname{Dom} \bolds\delta$; this does not follow directly from Proposition
\ref{pro:anticipativediff}, as the anticipative It\^o rule which was used
to obtain that result can yield integrands which are only in
$\mathbb{L}^{1,2}_{\rm loc}$.  We can verify directly, however, that
the integrand in this case is indeed in $\operatorname{Dom} \bolds\delta$, see Lemma
\ref{lem:skorisreal}.  Next, we must establish that $|U_{0,t}\nu|\,e_t^i$
is in $\mathbb{D}^{1,2}$ for every $i$.  Note that
$|U_{0,t}\nu|=\sum_i(U_{0,t}\nu)^i$, so $|U_{0,t}\nu|$ is in
$\mathbb{D}^\infty$.  Moreover, we establish in Lemma \ref{lem:regularpis}
that $e_t\in\mathbb{D}^{1,2}$ and that $\mall_re_t$ is a bounded
random variable for every $t$.  Hence it follows from Proposition
\ref{pro:chainrule} that $|U_{0,t}\nu|\,e_t^i\in\mathbb{D}^{1,2}$.
Consequently we can apply equation (\ref{eq:malladjoint}), and we obtain
\begin{eqnarray*}
&& \mathbf{E_Q}\biggl[|U_{0,t}\nu|\,e_t^*\int_0^t
D\tilde\pi_{r,t}(\pi_r)\cdot\Delta_H(\pi_r)\,dY_r\biggr]
\\
&&\qquad = \int_0^t\mathbf{E_Q}[(|U_{0,t}\nu|\,\mall_re_t^*+
\mall_r|U_{0,t}\nu|\,e_t^*) D\tilde\pi_{r,t}(\pi_r)\cdot\Delta_H(\pi_r)]\,dr
\\
&&\qquad = \int_0^t\mathbf{E_Q}[|U_{0,t}\nu|\,(\mall_r\pi_t-\mall_r\tilde\pi_t)^*
D\tilde\pi_{r,t}(\pi_r)\cdot\Delta_H(\pi_r)]\,dr
\\
&&\qquad\phantom{=}{} + \int_0^t\mathbf{E_Q}\Biggl[\sum_i(U_{r,t}HU_{0,r}\nu)^i\,e_t^*
D\tilde\pi_{r,t}(\pi_r)\cdot\Delta_H(\pi_r)\Biggr]\,dr.
\end{eqnarray*}
Now note that $|e_t^i|\le 1$, and that by Lemma \ref{lem:regularpis}
\[
|(\mall_r\pi_t-\mall_r\tilde\pi_t)^i|\le
|(\mall_r\pi_t)^i|+|(\mall_r\tilde\pi_t)^i|\le
\max_k|h^k| + \max_k|\tilde h^k|.
\]
Furthermore we can estimate
\[
\bigg|\frac{\sum_i(U_{r,t}HU_{0,r}\nu)^i}{|U_{0,t}\nu|}\bigg|\le
\frac{1}{|U_{0,t}\nu|}\sum_{i,j,k}U_{r,t}^{ij}\,|h^j|\,U_{0,r}^{jk}\nu^k
\le \max_k |h^k|,
\]
where we have used a.s. nonnegativity of the matrix elements of
$U_{0,r}$ and $U_{r,t}$ (this must be the case, as, for example, $U_{r,t}\mu$ has
nonnegative entries for any vector $\mu$ with nonnegative entries). Hence
we obtain
\begin{eqnarray*}
&& \mathbf{E_Q}\biggl[|U_{0,t}\nu|e_t^*\int_0^t
D\tilde\pi_{r,t}(\pi_r)\cdot\Delta_H(\pi_r)\,dY_r\biggr]
\\
&&\qquad \le \biggl(2\max_k|h^k|+\max_k|\tilde h^k|\biggr)
\int_0^t\mathbf{E_Q}|U_{0,t}\nu| |D\tilde\pi_{r,t}(\pi_r)\cdot\Delta_H(\pi_r)|\,dr.
\end{eqnarray*}
The result follows after straightforward manipulations.
\end{pf}

Unlike in the case $\tilde h=h$, we now have to deal also with second
derivatives of the filter with respect to its initial condition.
These can be estimated much in the same way as we dealt with the first
derivatives.

\begin{lem}
\label{lem:secondderivs}
Let $\tilde\lambda_{ij}>0$ $\forall i\ne j$ and $\mu\in\mathcal{S}^{d-1}$,
$v,w\in T\mathcal{S}^{d-1}$.  Then a.s.
\begin{eqnarray*}
&& |D^2\tilde\pi_{s,t}(\mu)\cdot v-D^2\tilde\pi_{s,t}(\mu)\cdot w|
\\
&&\qquad \le 2\sum_k\frac{|v^k+w^k|}{\mu^k}
\sum_j\frac{|v^j-w^j|}{\mu^j} \exp\biggl(-2(t-s)\min_{p,q\ne p}
\sqrt{\tilde\lambda_{pq}\tilde\lambda_{qp}}\biggr).
\end{eqnarray*}
Moreover, the result still holds if $\mu,v,w$ are
$\mathscr{F}_s^Y$-measurable random variables with values a.s. in
$\mathcal{S}^{d-1}$ and $T\mathcal{S}^{d-1}$, respectively.
\end{lem}

\begin{pf}
Proceeding as in the proof of Proposition \ref{pro:jacoboundsimple}, we
can calculate directly the second derivative of (\ref{eq:wronginirep}):
\[
\bigl(D^2\pi_t(\mu)\cdot v\bigr)^i =
-2\bigl(D\pi_t(\mu)\cdot v\bigr)^i \frac{\sum_j(v^j/\nu^j)\mathbf{P}(X_0=a_j|\mathscr{F}_t^Y)}
{\sum_j(\mu^j/\nu^j)\mathbf{P}(X_0=a_j|\mathscr{F}_t^Y)}.
\]
Setting $\mu=\nu$ and using the triangle inequality, we obtain
\[
|D^2\pi_t(\nu)\cdot v-D^2\pi_t(\nu)\cdot w|\le
2\sum_{i,j}\frac{|v^j(D\pi_t(\nu)\cdot v)^i
-w^j(D\pi_t(\nu)\cdot w)^i|}{\nu^j}.
\]
Another application of the triangle inequality and
using Proposition \ref{pro:jacoboundsimple} gives
\begin{eqnarray*}
&& |D^2\pi_t(\nu)\cdot v-D^2\pi_t(\nu)\cdot w|
\\
&&\qquad \le \sum_k\frac{|v^k+w^k|}{\nu^k}\,|D\pi_t(\nu)\cdot (v-w)|
+\sum_k\frac{|v^k-w^k|}{\nu^k}\,|D\pi_t(\nu)\cdot(v+w)|
\\
&&\qquad\le2\sum_k\frac{|v^k+w^k|}{\nu^k}\sum_j\frac{|v^j-w^j|}{\nu^j}
\exp\biggl(-2t\min_{p,q\ne p}\sqrt{\lambda_{pq}\lambda_{qp}}\biggr).
\end{eqnarray*}
We can now repeat the arguments of Corollary \ref{cor:forgetting} to
establish that the result still holds if we replace $\pi_{0,t}$ by
$\tilde\pi_{s,t}$, $\lambda_{pq}$ by $\tilde\lambda_{pq}$, and $\nu,v,w$ by
$\mathscr{F}_s^Y$-measurable random variables $\mu,v,w$.  This completes
the proof.
\end{pf}

We are now ready to complete the proof of Theorem \ref{thm:mainresult}.

\begin{pf*}{Proof of Theorem \ref{thm:mainresult}}
Set $\beta=2\min_{p,q\ne p}(\tilde\lambda_{pq}\tilde\lambda_{qp})^{1/2}$.
Let us collect all the necessary estimates.  First, we have
\[
\int_0^t \mathbf{E_P}
|D\tilde\pi_{r,t}(\pi_r)\cdot\Delta_\Lambda\pi_r|\,dr
\le\beta^{-1}\sup_{s\ge 0}\mathbf{E_P}\biggl(1/\min_k\pi_s^k\biggr)
|\Lambda^*-\tilde\Lambda^*|,
\]
as we showed in Section \ref{sec:forgetting}.  Next, we obtain
\[
\int_0^t\mathbf{E_P}
|D\tilde\pi_{r,t}(\pi_r)\cdot\Delta_H(\pi_r)|\,dr
\le\beta^{-1}\sup_{\pi\in\mathcal{S}^{d-1}}\sum_k|h^k-\tilde h^k+\tilde h^*\pi-h^*\pi|
\]
using Corollary \ref{cor:forgetting}.  Using the triangle inequality, we
can estimate this by
\[
\int_0^t\mathbf{E_P}|D\tilde\pi_{r,t}(\pi_r)\cdot\Delta_H(\pi_r)|\,dr
\le (d+1)\beta^{-1}|h-\tilde h|.
\]
Next, we estimate using Corollary \ref{cor:forgetting}
\begin{eqnarray*}
&& \int_0^t \mathbf{E_P}\big|D\tilde\pi_{r,t}(\pi_r)\cdot
\bigl(h^*\pi_r\,(H-h^*\pi_r)\pi_r-
\tilde h^*\pi_r\,(\tilde H-\tilde h^*\pi_r)\pi_r\bigr)\big|\,dr
\\
&&\qquad \le \beta^{-1}\sup_{\pi\in\mathcal{S}^{d-1}}
\sum_k|h^*\pi\,(h^k-h^*\pi)-\tilde h^*\pi(\tilde h^k-\tilde h^*\pi)|
\\
&&\qquad \le \beta^{-1}\biggl((d+1)\max_k|h^k|+
d\max_{k,\ell}|\tilde h^k-\tilde h^\ell|\biggr)|h-\tilde h|,
\end{eqnarray*}
where we have used the estimate
\begin{eqnarray*}
&&\sum_k|h^*\pi\,(h^k-h^*\pi)-\tilde h^*\pi\,(\tilde h^k-\tilde h^*\pi)|
\\
&&\qquad \le |h^*\pi| \sum_k|h^k-\tilde h^k+\tilde h^*\pi-h^*\pi|
+|h^*\pi-\tilde h^*\pi|\sum_k|\tilde h^k-\tilde h^*\pi|
\\
&&\qquad \le (d+1)\max_k|h^k|~|h-\tilde h|+
|h-\tilde h|\sum_k|\tilde h^k-\tilde h^*\pi|
\\
&&\qquad \le \biggl((d+1)\max_k|h^k|+ d\max_{k,\ell}|\tilde h^k-\tilde h^\ell|\biggr)|h-\tilde h|.
\end{eqnarray*}
Next we estimate using Lemma \ref{lem:secondderivs}
\begin{eqnarray*}
&& \tfrac{1}{2}\int_0^t \mathbf{E_P}|
D^2\tilde\pi_{r,t}(\pi_r)\cdot(H-h^*\pi_r)\pi_r
-D^2\tilde\pi_{r,t}(\pi_r)\cdot(\tilde H-\tilde h^*\pi_r)\pi_r|\,dr
\\
&&\qquad \le \beta^{-1}\sup_{\pi\in\mathcal{S}^{d-1}}
\sum_k|h^k-h^*\pi+\tilde h^k-\tilde h^*\pi|
\sum_j|h^j-\tilde h^j+\tilde h^*\pi-h^*\pi|
\\
&&\qquad \le  d(d+1)\beta^{-1}\biggl(\max_{k,\ell}|h^k-h^\ell|+
\max_{k,\ell}|\tilde h^k-\tilde h^\ell|\biggr)|h-\tilde h|.
\end{eqnarray*}
We have now estimated all the terms in Lemma \ref{lem:filtererror}, and
hence we have bounded $\mathbf{E_P}\|e_t\|^2=
\mathbf{E_P}\|\pi_t(\nu)-\tilde\pi_t(\nu)\|^2$.  It remains to allow for
misspecified initial conditions.  To this end, we estimate
\begin{eqnarray*}
&&\|\pi_t(\nu)-\tilde\pi_t(\mu)\|^2
\\
&&\qquad \le
\|e_t\|^2+\|\tilde\pi_t(\nu)-\tilde\pi_t(\mu)\|
\bigl(\|\tilde\pi_t(\nu)-\tilde\pi_t(\mu)\|+2\|\pi_t(\nu)-\tilde\pi_t(\nu)\|\bigr).
\end{eqnarray*}
Hence we obtain using the equivalence of finite-dimensional norms
$\|x\|\le K_{21}\,|x|$
\[
\|\pi_t(\nu)-\tilde\pi_t(\mu)\|^2\le
\|e_t\|^2+ 6K_{21}\,|\tilde\pi_t(\nu)-\tilde\pi_t(\mu)|
\]
where we have used that the simplex is contained in the
$(d-1)$-dimensional unit sphere, so $\|\mu_1-\mu_2\|\le 2$ $\forall
\mu_1,\mu_2\in\Delta^{d-1}$.  The statement of the theorem now follows
directly from Lemma \ref{lem:filtererror}, Proposition
\ref{pro:forgetting} and the estimates above.
\end{pf*}

\appendix

\section{\texorpdfstring{Anticipative stochastic calculus}{Appendix A: Anticipative stochastic calculus}}
\label{sec:malliavin}

The goal of this appendix is to recall briefly the main results of the
Malliavin calculus, Skorokhod integrals and anticipative stochastic
calculus that are needed in the proofs.  In our application of the theory
we wish to deal with functionals of the observation process
$(Y_t)_{t\in[0,T]}$, where $T$ is some finite time (usually we will
calculate integrals from $0$ to $t$, so we can choose any $T>t$).  Recall
from Section~\ref{sec:prelims} that $Y$ is an $\mathscr{F}_t^Y$-Wiener
process under the measure $\mathbf{Q}$; it will thus be convenient to work
always under $\mathbf{Q}$, as this puts us directly in the framework
used, for example, in \cite{nualart}.  As the theory described below is defined
$\mathbf{Q}$-a.s. and as $\mathbf{P}\sim\mathbf{Q}$, the corresponding
properties under $\mathbf{P}$ are unambiguously obtained by using equation
(\ref{eq:refmeasure}).  We will presume this setup whenever the theory
described here is applied.

A smooth random variable $F$ is one of the form
$f(Y(h_1),\ldots,Y(h_n))$, where $Y(h)$ denotes the Wiener integral of the
deterministic function $h\in L^2([0,T])$ with respect to $Y$ and $f$ is a
smooth function which is of polynomial growth together with all its
derivatives.  For smooth $F$ the Malliavin derivative $\mall F$ is
defined by
\[
\mall_tF=\sum_{i=1}^n\frac{\partial f}{\partial x^i}
(Y(h_1),\ldots,Y(h_n))h_i(t).
\]
The Malliavin derivative $\mall$ can be shown \cite{nualart}, page 26,
to be closeable as an operator from
$L^p(\Omega,\mathscr{F}_T^Y,\mathbf{Q})$ to
$L^p(\Omega,\mathscr{F}_T^Y,\mathbf{Q};L^2([0,T]))$ for any $p\ge 1$, and
we denote the domain of $\mall$ in $L^p(\Omega)$ by $\mathbb{D}^{1,p}$
[for notational convenience we will drop the measure $\mathbf{Q}$
and $\sigma$-algebra $\mathscr{F}_T^Y$ throughout this section,
where it is understood that $L^p(\Omega)$ denotes
$L^p(\Omega,\mathscr{F}_T^Y,\mathbf{Q})$, etc.].
More generally, we consider iterated derivatives
$\mall^kF\in L^p(\Omega;L^2([0,T]^k))$ defined by
$\mall^k_{t_1,\ldots,t_k}F=\mall_{t_1}\cdots\mall_{t_k}F$, and the domain
of $\mall^k$ in $L^p(\Omega)$ is denoted by $\mathbb{D}^{k,p}$.
The domains $\mathbb{D}^{k,p}$ can also be localized
(\cite{nualart}, pages 44--45), and we denote the corresponding localized domains by
$\mathbb{D}^{k,p}_{\rm loc}$.  Finally, we define the useful class
$\mathbb{D}^\infty=\bigcap_{p\ge 1}\bigcap_{k\ge 1}\mathbb{D}^{k,p}$.

We will use two versions of the chain rule for the Malliavin derivative.

\begin{prop}
\label{pro:chainrule}
Let $\varphi\dvtx\mathbb{R}^m\to\mathbb{R}$ be $C^1$ and
$F=(F^1,\ldots,F^m)$ be a random vector with components in
$\mathbb{D}^{1,2}$.  Then $\varphi(F)\in\mathbb{D}^{1,2}_{\rm loc}$ and
\[
\mall\varphi(F)=\sum_{i=1}^m\frac{\partial\varphi}{\partial x^i}
(F)\,\mall F^i.
\]
If $\varphi(F)\in L^2(\Omega)$ and $\mall\varphi(F)\in L^2(\Omega\times [0,T])$,
then $\varphi(F)\in\mathbb{D}^{1,2}$.  These results
still hold if $F$ a.s. takes values in an open domain
$V\subset\mathbb{R}^m$ and $\varphi$ is $C^1(V)$.
\end{prop}

The first (local) statement is \cite{nualartpardoux}, Proposition 2.9; the
second statement can be proved in the same way as
\cite{oconekaratzas}, Lemma A.1, and the proofs are easily adapted to the case where
$F$ a.s. takes values in some domain.  The next result is from
\cite{nualart}, page 62:

\begin{prop}
\label{pro:chainrulesmooth}
Let $\varphi\dvtx\mathbb{R}^m\to\mathbb{R}$ be a smooth function which is of
polynomial growth together with all its derivatives, and let
$F=(F^1,\ldots,F^m)$ be a random vector with components in
$\mathbb{D}^\infty$.  Then $\varphi(F)\in\mathbb{D}^\infty$ and the usual
chain rule holds.  This implies that $\mathbb{D}^\infty$ is an algebra,
that is, $FG\in\mathbb{D}^\infty$ for $F,G\in\mathbb{D}^\infty$.
\end{prop}

The following result follows from \cite{nualart}, page 32
(here $[s,t]^c=[0,T]\backslash [s,t]$).

\begin{lem}
\label{lem:adapted}
If $F\in\mathbb{D}^{1,2}$ is $\mathscr{F}^Y_{[s,t]}$-measurable, then
$\mall F=0$ a.e. in $\Omega\times[s,t]^c$.
\end{lem}

It is useful to be able to calculate explicitly the Malliavin derivative
of the solution of a stochastic differential equation.  Consider
$dx_t = f(x_t)\,dt+\sigma(x_t)\,dY_t$, $x_0\in\mathbb{R}^m$,
where $f(x)$ and $\sigma(x)$ are smooth functions of $x$ with bounded
derivatives of all orders.  It is well known that such equations generate
a smooth stochastic flow of diffeomorphisms $x_t=\xi_t(x)$ \cite{kunita}.
We now have the following result.

\begin{prop}
\label{pro:sdemalliavin}
All components of $x_t$ belong to $\mathbb{D}^\infty$ for every
$t\in[0,T]$. We have $\mall_rx_t=D\xi_t(x_0)D\xi_r(x_0)^{-1}\sigma(x_r)$ a.e. $r<t$, where
$(D\xi_t(x))^{ij}=\partial\xi_t^i(x)/\partial x^j$ is the Jacobian matrix
of the flow, and $\mall_rx_t=0$ a.e. $r>t$.
\end{prop}

The first statement is given in \cite{nualart}, Theorem 2.2.2, page 105,
the second on \cite{nualart}, equation (2.38), page 109, the third
follows from adaptedness (Lemma \ref{lem:adapted}).

We now consider $\mall$ as a closed operator from
$L^2(\Omega)$ to $L^2(\Omega\times[0,T])$ with domain $\mathbb{D}^{1,2}$.
Its Hilbert space adjoint $\bolds\delta=\mall^*$ is well defined in the usual
sense as a closed operator from $L^2(\Omega\times[0,T])$ to $L^2(\Omega)$,
and we denote its domain by $\operatorname{Dom}\bolds\delta$.  The operator $\bolds\delta$ is
called the Skorokhod integral, and coincides with the It\^o integral on
the subspace $L^2_a(\Omega\times[0,T])\subset \operatorname{Dom}\bolds\delta$ of
adapted square integrable processes (\cite{nualart}, Proposition 1.3.4, page 41).
$\bolds\delta$ is thus an extension of the It\^o integral to a class of possibly
anticipative integrands.  To emphasize this point we will write
\[
\bolds\delta \bigl(uI_{[s,t]}\bigr)=\int_s^tu_r\,dY_r,\qquad
uI_{[s,t]}\in \operatorname{Dom}\bolds\delta.
\]
The Skorokhod integral has the following properties.  First, its
expectation vanishes $\mathbf{E_Q}\bolds\delta (u)=0$ if $u\in\operatorname{Dom}\bolds\delta$.
Second, by its definition as the adjoint of $\mall$ we have
\begin{equation}
\label{eq:malladjoint}
\mathbf{E_Q}(F\bolds\delta(u))=\mathbf{E_Q}\biggl[
\int_0^T(\mall_tF)u_t\,dt\biggr]
\end{equation}
if $u\in\operatorname{Dom}\bolds\delta$, $F\in\mathbb{D}^{1,2}$. We will also use the
following result, the proof of which proceeds in exactly the same way as
its one-dimensional counterpart (\cite{nualart}, page 40).

\begin{lem}
\label{lem:bringintoskor}
If $u$ is an $n$-vector of processes in $\operatorname{Dom}\bolds\delta$ and
$F$ is an \mbox{$m\times n$-matrix} of random variables in
$\mathbb{D}^{1,2}$ such that $\mathbf{E_Q}\int_0^T\|Fu_t\|^2\,dt<\infty$, then
\[
\int_0^TFu_t\,dY_t=F\int_0^Tu_t\,dY_t-\int_0^T
(\mall_tF)u_t\,dt
\]
in the sense that $Fu\in\operatorname{Dom}\bolds\delta$ iff the right-hand side of this
expression is in $L^2(\Omega)$.
\end{lem}

As it is difficult to obtain general statements for integrands in
$\operatorname{Dom}\bolds\delta$, it is useful to single out restricted classes of integrands
that are easier to deal with.  To this end, define the spaces
$\mathbb{L}^{k,p}=L^p([0,T];\mathbb{D}^{k,p})$ for $k\ge 1$, $p\ge 2$.
Note that $\mathbb{L}^{k,p}\subset\mathbb{L}^{1,2}\subset \operatorname{Dom}\bolds\delta$
\cite{nualart}, page 38.  Moreover, the domains $\mathbb{L}^{k,p}$ can be
localized to $\mathbb{L}^{k,p}_{\rm loc}$ (\cite{nualart}, pages 43--45).
We can now state an It\^o change of variables formula for Skorokhod
integrals, see \cite{nualart,nualartpardoux,oconepardoux}.  The extension
to processes that a.s. take values in some domain is straightforward
through localization.

\begin{prop}
\label{pro:anticipatingito}
Consider an $m$-dimensional process of the form
\[
x_t=x_0+\int_0^t v_s\,ds+\int_0^tu_s\,dY_s,
\]
where we assume that $x_t$ has a continuous version and
$x_0\in(\mathbb{D}^{1,4}_{\rm loc})^m$, $v\in(\mathbb{L}^{1,4}_{\rm
loc})^m$, and $u\in(\mathbb{L}^{2,4}_{\rm loc})^m$.  Let
$\varphi:\mathbb{R}^m\to\mathbb{R}$ be a $C^2$ function.  Then
\begin{eqnarray*}
\varphi(x_t) &=& \varphi(x_0)+\int_0^t D\varphi(x_s)v_s\,ds+
\int_0^t D\varphi(x_s)u_s\,dY_s
\\
&&{} +\tfrac{1}{2} \int_0^t (D^2\varphi(x_s)\nab_sx_s,u_s)\,ds,
\end{eqnarray*}
where $\nab_sx_s=\lim_{\varepsilon\searrow
0}\mall_s(x_{s+\varepsilon}+x_{s-\varepsilon})$,
$D\varphi(x_s)u_s=\sum_i(\partial\varphi/\partial
x^i)(x_s)u_s^i$,\break
$(D^2\varphi(x_s)\nab_sx_s,u_s)=\sum_{ij}(\partial^2\varphi/\partial
x^i\,\partial x^j)(x_s)u_s^i\nab_sx_s^j$.  The result still holds
if $x_s$ a.s. takes values in an open domain $V\subset\mathbb{R}^m$
$\forall s\in[0,t]$ and $\varphi$ is $C^2(V)$.
\end{prop}

\section{\texorpdfstring{Some technical results}{Appendix B: Some technical results}}
\label{sec:proofs}

\begin{lem}
\label{lem:msqintg}
The following equality holds:
\begin{eqnarray*}
&& \tilde U_{0,t}\int_0^s\tilde U_{0,r}^{-1}(H-\tilde H)U_{0,r}\nu\,dY_r
\\
&&\qquad=\int_0^s\tilde U_{r,t}(H-\tilde H)U_{0,r}\nu\,dY_r+
\int_0^s \tilde U_{r,t}\tilde H(H-\tilde H)U_{0,r}\nu\,dr.
\end{eqnarray*}
The integral on the left-hand side is an It\^o integral, on the
right-hand side a Skorokhod integral.
\end{lem}

\begin{pf}
We have already established in the proof of Proposition
\ref{pro:anticipativediff} that the matrix elements of $\tilde U_{0,t}$ are
in $\mathbb{D}^\infty\subset\mathbb{D}^{1,2}$.  Moreover,
\begin{eqnarray*}
&& \mathbf{E_Q}\|\tilde U_{r,t}(H-\tilde H)U_{0,r}\nu\|^2
\\
&&\qquad\le \|H-\tilde H\|^2\,\mathbf{E_Q}(\|\tilde U_{r,t}\|^2\,\|U_{0,r}\|^2)
\\
&&\qquad\le \|H-\tilde H\|^2\,\sqrt{\mathbf{E_Q}\|\tilde U_{r,t}\|^4\,
\mathbf{E_Q}\|U_{0,r}\|^4}
\\
&&\qquad\le C_4^4\,\|H-\tilde H\|^2\,\sqrt{\mathbf{E_Q}\|\!|\tilde U_{r,t}\|\!|^4_4\,
\mathbf{E_Q}\|\!|U_{0,r}\|\!|^4_4},
\end{eqnarray*}
where we have used the Cauchy--Schwarz inequality and $\|\nu\|\le 1$ for
$\nu\in\mathcal{S}^{d-1}$.  Here $\|\!|U\|\!|_p=(\sum_{ij}U_{ij}^p)^{1/p}$ is
the elementwise $p$-norm of $U$, $\|U\|$ is the usual matrix $2$-norm,
and $C_p$ matches the norms $\|U\|\le C_p\|\!|U\|\!|_p$ (recall
that all norms on a finite-dimensional space are equivalent).  As
$U_{0,r},\tilde U_{r,t}$ are solutions of linear stochastic differential
equations, standard estimates give for any integer $p\ge 2$
\begin{eqnarray*}
\mathbf{E_Q}\biggl(\sup_{0\le r\le t}\|\!|\tilde U_{r,t}\|\!|_p^p\biggr)
&\le & D_1(p)<\infty,
\\
\mathbf{E_Q}\biggl(\sup_{0\le r\le t}\|\!|U_{0,r}\|\!|_p^p\biggr)
&\le& D_2(p)<\infty,
\end{eqnarray*}
and we obtain
\[
\int_0^s\mathbf{E_Q}\|\tilde U_{r,t}(H-\tilde H)U_{0,r}\nu\|^2\,dr\le
s\sup_{0\le r\le s}\mathbf{E_Q}\|\tilde U_{r,t}(H-\tilde H)U_{0,r}\nu\|^2<\infty.
\]
Hence we can apply Lemma \ref{lem:bringintoskor} to obtain the
result.  By a similar calculation we can establish that the right-hand
side of the expression in Lemma \ref{lem:bringintoskor} for our case is
square integrable, so that the Skorokhod integral is well defined.
\end{pf}

\begin{lem}
\label{lem:canapplyito}
The anticipating It\^o rule with $\Sigma(x)=x/|x|$ can be applied to
\[
\tilde U_{s,t}U_{0,s}\nu=\tilde U_{0,t}\nu
+\int_0^s\tilde U_{r,t}(\Lambda^*-\tilde\Lambda^*)U_{0,r}\nu\,dr
+\int_0^s\tilde U_{r,t}(H-\tilde H)U_{0,r}\nu\,dY_r.
\]
\end{lem}

\begin{pf}
Clearly the Skorokhod integral term has a.s. continuous sample paths, as
both $\tilde U_{s,t}U_{0,s}\nu$ and the time integrals do; moreover,
$\tilde U_{0,t}\nu\in(\mathbb{D}^\infty)^d$.  In order to be able to apply
Proposition \ref{pro:anticipatingito}, it remains to check the
technical conditions $v_r=\tilde U_{r,t}(\Lambda^*-\tilde\Lambda^*)U_{0,r}\nu\in
(\mathbb{L}^{1,4})^d$, $u_r=\tilde U_{r,t}(H-\tilde H)U_{0,r}\nu\in
(\mathbb{L}^{2,4})^d$.

As $\mathbb{D}^\infty$ is an algebra, $u_t$ and $v_t$ take values in
$\mathbb{D}^\infty$.  Moreover, we can establish exactly as in the proof
of Lemma \ref{lem:msqintg} that $u$ and $v$ are in $L^4(\Omega\times [0,t])$.
To complete the proof we must establish that
\begin{eqnarray*}
&& \sum_i\int_0^t\mathbf{E_Q}\biggl[\int_0^t
(\mall_su_r^i)^2\,ds\biggr]^{2}\,dr<\infty,
\\
&& \sum_i\int_0^t\mathbf{E_Q}\biggl[\int_0^t
(\mall_sv_r^i)^2\,ds \biggr]^{2}\,dr<\infty,
\end{eqnarray*}
thus ensuring that $u,v\in(\mathbb{L}^{1,4})^d$, and
\[
\sum_i\int_0^t\mathbf{E_Q}\biggl[\int_0^t\int_0^t
(\mall_\sigma\mall_su_r^i)^2\,ds\,d\sigma\biggr]^{2}\,dr<\infty
\]
which ensures that $u\in(\mathbb{L}^{2,4})^d$.  Using the Cauchy--Schwarz
inequality we have
\begin{eqnarray*}
&& \sum_i\int_0^t\mathbf{E_Q}\biggl[\int_0^t
(\mall_su_r^i)^2\,ds\biggr]^{2}\,dr
\\
&&\qquad \le t\int_0^t\int_0^t\mathbf{E_Q}\|\mall_su_r\|_4^4\,ds\,dr
\le t^3\sup_{0\le r,s\le t}\mathbf{E_Q}\|\mall_su_r\|_4^4,
\end{eqnarray*}
and similarly for $v$.  Moreover, we obtain
\[
\sum_i\int_0^t\mathbf{E_Q}\biggl[\int_0^t\int_0^t
(\mall_\sigma\mall_su_r^i)^2\,ds\,d\sigma\biggr]^{2}\,dr\le
t^5\sup_{0\le r,s,\sigma\le t}
\mathbf{E_Q}\|\mall_\sigma\mall_su_r\|_4^4.
\]
But using the chain rule Proposition \ref{pro:chainrulesmooth} we can
easily establish that
\[
\mall_su_r=\cases{
\tilde U_{r,t}(H-\tilde H)U_{s,r}HU_{0,s}\nu, &\quad \mbox{a.e.} $0<s<r<t$,
\cr
\tilde U_{s,t}\tilde H\tilde U_{r,s}(H-\tilde H)U_{0,r}\nu, &\quad\mbox{a.e.} $0<r<s<t$,}
\]
and similarly
\[
\mall_\sigma\mall_su_r=\cases{
\tilde U_{r,t}(H-\tilde H)U_{s,r}HU_{\sigma,s}HU_{0,\sigma}\nu, &
\quad\mbox{a.e.} $0<\sigma<s<r<t$,
\cr
\tilde U_{r,t}(H-\tilde H)U_{\sigma,r}HU_{s,\sigma}HU_{0,s}\nu, &
\quad \mbox{a.e.} $0<s<\sigma<r<t$,
\cr
\tilde U_{\sigma,t}\tilde H\tilde U_{r,\sigma}(H-\tilde H)U_{s,r}HU_{0,s}\nu, &
\quad\mbox{a.e.} $0<s<r<\sigma<t$,
\cr
\tilde U_{s,t}\tilde H\tilde U_{r,s}(H-\tilde H)U_{\sigma,r}HU_{0,\sigma}\nu, &
\quad \mbox{a.e.} $0<\sigma<r<s<t$,
\cr
\tilde U_{s,t}\tilde H\tilde U_{\sigma,s}\tilde H\tilde U_{r,\sigma}(H-\tilde H)U_{0,r}\nu, &
\quad \mbox{a.e.} $0<r<\sigma<s<t$,
\cr
\tilde U_{\sigma,t}\tilde H\tilde U_{s,\sigma}\tilde H\tilde U_{r,s}(H-\tilde H)U_{0,r}\nu, &
\quad\mbox{a.e.} $0<r<s<\sigma<t$.}
\]
The desired estimates now follow as in the proof of Lemma \ref{lem:msqintg}.
\end{pf}

\begin{lem}
\label{lem:skorisreal}
The Skorokhod integrand obtained by applying the anticipative It\^o
formula as in Lemma \emph{\ref{lem:canapplyito}} is in $\operatorname{Dom} \bolds\delta$.
\end{lem}

\begin{pf}
We use the notation $\rho_r=U_{0,r}\nu$.  The Skorokhod integral in
question is
\[
\int_0^s D\Sigma(\tilde U_{r,t}\rho_r)\tilde U_{r,t}(H-\tilde H)\rho_r\,dY_r=
\int_0^s f_r\,dY_r.
\]
To establish $f\in\operatorname{Dom} \bolds\delta$, it suffices to show that
$f\in\mathbb{L}^{1,2}$.  We begin by showing
\begin{eqnarray*}
&& |D\Sigma(\tilde U_{r,t}\rho_r)\tilde U_{r,t}(H-\tilde H)\rho_r|
\\
&&
\qquad =\sum_i\Bigg|\sum_{j,k}\frac{\delta^{ij}-
\Sigma^i(\tilde U_{r,t}\rho_r)}{|\tilde U_{r,t}\rho_r|}
\tilde U_{r,t}^{jk}(h^k-\tilde h^k)\rho_r^k\Bigg|
\\
&&\qquad\le \frac{1}{|\tilde U_{r,t}\rho_r|}\sum_{i,j,k}
\tilde U_{r,t}^{jk}\,|h^k-\tilde h^k|\,\rho_r^k
\\
&&\qquad\le \frac{\max_k|h^k-\tilde h^k|}{|\tilde U_{r,t}\rho_r|}\sum_{i,j,k}
\tilde U_{r,t}^{jk}\rho_r^k
\\
&&\qquad =d\max_k|h^k-\tilde h^k|,
\end{eqnarray*}
where we have used the triangle inequality,
$|\delta^{ij}-\Sigma^{i}(x)|\le 1$ for any $x\in\mathbb{R}^d_{++}$, and
the fact that $U_{r,t}$ and $\rho_r$ have nonnegative entries a.s.
Hence $f_r$ is a bounded process.  Similarly, we will show that
$\mall_sf_r$ is a bounded process.  Note that $f_r$ is a smooth function
on $\mathbb{R}^d_{++}$ of positive random variables in
$\mathbb{D}^\infty$; hence we can apply the chain rule Proposition
\ref{pro:chainrule}.  This gives
\[
(\mall_sf_r)^i=\cases{\displaystyle
\sum_{jk} D^2\Sigma^{ijk}(\tilde U_{r,t}\rho_r)
\bigl(\tilde U_{r,t}(H-\tilde H)\rho_r\bigr)^j(\tilde U_{r,t}U_{s,r}H\rho_s)^k
\cr\displaystyle
\qquad{} + \sum_j D\Sigma^{ij}(\tilde U_{r,t}\rho_r)
\bigl(\tilde U_{r,t}(H-\tilde H)U_{s,r}H\rho_s\bigr)^j, & \quad \mbox{a.e.} $s<r$,
\cr\displaystyle
\sum_{jk} D^2\Sigma^{ijk}(\tilde U_{r,t}\rho_r)
\bigl(\tilde U_{r,t}(H-\tilde H)\rho_r\bigr)^j (\tilde U_{s,t}\tilde H\tilde U_{r,s}\rho_r)^k
\cr\displaystyle
\qquad{} + \sum_j D\Sigma^{ij}(\tilde U_{r,t}\rho_r)
\bigl(\tilde U_{s,t}\tilde H\tilde U_{r,s}(H-\tilde H)\rho_r\bigr)^j, &
\quad\mbox{a.e.} $s>r$.}
\]
Proceeding exactly as before, we find that
$\mall f\in L^\infty(\Omega\times [0,t]^2)$.  But then by Proposition
\ref{pro:chainrule} we can conclude that $\mall_sf_r\in\mathbb{D}^{1,2}$
for a.e.\ $(s,t)\in[0,t]^2$,  and in particular $f\in\mathbb{L}^{1,2}$.
Hence the proof is complete.
\end{pf}

\begin{lem}
\label{lem:regularpis}
$\mall_r\pi_s=D\pi_{r,s}(\pi_r)\cdot (H-h^*\pi_r)\pi_r$ a.e. $r<s$,
$\mall_r\pi_s=0$ a.e. $r>s$.  Moreover $|(\mall_r\pi_s)^i|\le
\max_k|h^k|$ for every $i$. The equivalent results hold for
$\mall_r\tilde\pi_s$.  In particular, this implies that $\pi_s$ and $\tilde\pi_s$
are in $\mathbb{D}^{1,2}$.
\end{lem}

\begin{pf}
The case $r>s$ is immediate from adaptedness of $\pi_s$.  For $r<s$,
apply the chain rule to $\pi_s=\Sigma(U_{0,s}\nu)\in\mathbb{D}^{1,2}_{\rm
loc}$.  Boundedness of the resulting expression follows, for example, as in the
proof of Lemma \ref{lem:skorisreal}, and hence it follows that
$\pi_s\in\mathbb{D}^{1,2}$.
\end{pf}

\subsection*{Acknowledgment}

R. van Handel thanks
P.~S.~Krishnaprasad of the University of Maryland for hosting a visit to
the Institute for Systems Research, during which the this work was initiated.

\printaddresses

\begin{thebibliography}{10}

\bibitem{atar}
\textsc{Atar, R.} and \textsc{Zeitouni, O.} (1997).
Lyapunov exponents for finite state nonlinear filtering.
\textit{SIAM J. Control Optim.} \textbf{35} 36--55.
\MR{1430282}

\bibitem{BaxChiLip}
\textsc{Baxendale, P., Chigansky, P.} and \textsc{Liptser, R.} (2004).
Asymptotic stability of the Wonham filter: Ergodic and nonergodic signals.
\textit{SIAM J. Control Optim.} \textbf{43} 643--669.
\MR{2086177}

\bibitem{BhatKalKar1}
\textsc{Bhatt, A.~G., Kallianpur, G.} and \textsc{Karandikar, R.~L.} (1995).
Uniqueness and robustness of solution of measure-valued equations of
nonlinear filtering. \textit{Ann. Probab.} \textbf{23} 1895--1938.
\MR{1379173}

\bibitem{BhatKalKar2}
\textsc{Bhatt, A.~G., Kallianpur, G.} and \textsc{Karandikar, R.~L.} (1999).
Robustness of the nonlinear filter.
\textit{Stochastic Process. Appl.} \textbf{81} 247--254.
\MR{1694557}

\bibitem{brigo}
\textsc{Brigo, D., Hanzon, B.} and \textsc{Le~Gland, F.} (1999).
Approximate nonlinear filtering by projection on exponential
manifolds of densities. \textit{Bernoulli} \textbf{5} 495--534.
\MR{1693600}

\bibitem{kushner1}
\textsc{Budhiraja, A.} and \textsc{Kushner, H.~J.} (1998).
Robustness of nonlinear filters over the infinite time interval.
\textit{SIAM J. Control Optim.} \textbf{36} 1618--1637.
\MR{1626884}

\bibitem{kushner2}
\textsc{Budhiraja, A.} and \textsc{Kushner, H.~J.} (1999).
Approximation and limit results for nonlinear filters over an
infinite time interval. \textit{SIAM J. Control Optim.} \textbf{37} 1946--1979.
\MR{1720146}

\bibitem{daprato}
\textsc{Da~Prato, G., Fuhrman, M.} and \textsc{Malliavin, P.} (1999).
Asymptotic ergodicity of the process of conditional law in some
problem of non-linear filtering. \textit{J. Funct. Anal.} \textbf{164} 356--377.
\MR{1695555}

\bibitem{elliott}
\textsc{Elliott, R.~J., Aggoun, L.} and \textsc{Moore, J.~B.} (1995).
\textit{Hidden Markov Models}. Springer, New York.
\MR{1323178}

\bibitem{GuoYin}
\textsc{Guo, X.} and \textsc{Yin, G.} (2006).
The Wonham filter with random parameters: Rate of convergence and error bounds.
\textit{IEEE Trans. Automat. Control} \textbf{51} 460--464.
\MR{2205683}

\bibitem{kunita}
\textsc{Kunita, H.} (1984).
Stochastic differential equations and stochastic flows of diffeomorphisms.
\textit{\'Ecole d'\'Et\'e de Probabilit\'es de Saint-Flour
XII---1982. Lecture Notes in Math.} \textbf{1097} 143--303.
Springer, Berlin.
\MR{0876080}

\bibitem{legland1}
\textsc{Le~Gland, F.} and \textsc{Oudjane, N.} (2003).
A robustification approach to stability and to uniform particle
approximation of nonlinear filters: The example of pseudo-mixing signals.
\textit{Stochastic Process. Appl.} \textbf{106} 279--316.
\MR{1989630}

\bibitem{legland2}
\textsc{Le~Gland, F.} and \textsc{Oudjane, N.} (2004).
Stability and uniform approximation of nonlinear filters using the
Hilbert metric and application to particle filters.
\textit{Ann. Appl. Probab.} \textbf{14} 144--187.
\MR{2023019}

\bibitem{liptser}
\textsc{Liptser, R.~S.} and \textsc{Shiryaev, A.~N.} (2001).
\textit{Statistics of Random Processes} I. Springer, Berlin.
\MR{1800857}

\bibitem{nualart}
\textsc{Nualart, D.} (1988).
\textit{The Malliavin Calculus and Related Topics}.
Springer, New York.
\MR{2200233}

\bibitem{nualartpardoux}
\textsc{Nualart, D.} and \textsc{Pardoux, \'E.} (1995).
Stochastic calculus with anticipating integrands.
\textit{Probab. Theory Related Fields} \textbf{78} 535--581.
\MR{0950346}

\bibitem{oconekaratzas}
\textsc{Ocone, D.~L.} and \textsc{Karatzas, I.} (1991).
A generalized Clark representation formula, with application to
optimal portfolios. \textit{Stochastics Stochastics Rep.} \textbf{34} 187--220.
\MR{1124835}

\bibitem{oconepardoux}
\textsc{Ocone, D.} and \textsc{Pardoux, \'E.} (1989).
A generalized It\^{o}--Ventzell formula. Application to a class of
anticipating stochastic differential equations.
\textit{Ann. Inst. H. Poincar\'e Probab. Statist.} \textbf{25} 39--71.
\MR{0995291}

\bibitem{papa}
\textsc{Papavasiliou, A.} (2006).
Parameter estimation and asymptotic stability in stochastic filtering.
\textit{Stochastic Process. Appl.} \textbf{116} 1048--1065.
\MR{2238613}

\bibitem{protter}
\textsc{Protter, P.~E.} (2004).
\textit{Stochastic Integration and Differential Equations}, 2nd ed.
Springer, Berlin.
\MR{2020294}

\bibitem{rogersw}
\textsc{Rogers, L.~C.~G.} and \textsc{Williams, D.} (2000).
\textit{Diffusions, Markov Processes, and Martingales}. \textbf{2}.
\textit{It\^{o} Calculus}. Cambridge Univ. Press.
\MR{1780932}

\bibitem{wonham}
\textsc{Wonham, W.~M.} (1965).
Some applications of stochastic differential equations to optimal
nonlinear filtering. \textit{J. Soc. Indust. Appl. Math. Ser. A Control}
\textbf{2} 347--369.
\MR{0186472}

\end{thebibliography}
\end{document}